\documentclass[preprint,12pt]{elsarticle}
\usepackage[margin=1in]{geometry}
\usepackage{setspace}
\usepackage{adjustbox}
\usepackage{array}
\usepackage{mathrsfs}
\usepackage[normalem]{ulem}
\usepackage{dirtytalk}
\usepackage{natbib}
\setcitestyle{authoryear,open={(},close={)}}
\usepackage{pgfplots}
\usepackage{multirow}
\usepackage{lscape} 
\usepackage{arydshln}
\usepackage{tikz}
\usetikzlibrary{arrows.meta}
\usepackage{pifont}
\usepackage{enumitem}

\usepackage{color, colortbl}
\usepackage{pgfplots}
\usepackage{pgfplotstable}
\definecolor{Gray}{gray}{0.9}
\usepackage{float}
\usetikzlibrary[positioning]
\usetikzlibrary{patterns}
\usepackage[toc,page]{appendix}
\usepackage{tkz-euclide}
\usepackage{algorithm,amsmath,lipsum,xcolor,caption}
\usepackage{graphicx}
\usepackage{caption}
\usepackage{subcaption}
\usepackage{setspace}
\usepackage{forest}
\usepackage{xcolor,colortbl}
\usepackage{arydshln}
\usepackage{amssymb}

\newtheorem{theorem}{Theorem}
\newtheorem{lemma}{Lemma}
\newtheorem{definition}{Definition}

\newtheorem{remark}{Remark}

\newtheorem{example}{Example}

\begin{document}

\begin{frontmatter}

\title{Stable allocations for choice-based collaborative price setting}
% Other titles:
%\title{Stable allocations for collaborative price setting initiatives: \\ \normalsize{merging cooperative game theory and discrete-choice theory}}
%\color{red}A cooperative game theory approach for choice-based price setting \color{black}
%Collaboration in a choice-based pricing \textcolor{red}{setting}:\\ a cooperative game approach
%Stable allocations for collaborative initiatives in choice-based price settings
%Stable allocation for cooperation problems in choice-based price settings
% Core Allocations for Cooperation Pricing Problems \\ under Multinomial Logit
% Stable collaboration amongst start-up companies: combining cooperative game and discrete choice theories
%Stable allocations for cooperation problems in choice-based price settings.
% Unlocking benefits of micromobility via collaboration: combining cooperative game and discrete choice theories
%Supporting startups operators via collaboration:
%A Cooperative Game with Discrete Choice Modelling}
%How to allocate costs when transport companies collaborate? A cooperative game theory approach
%Cooperartive Pricing under Multinomial Logit Demand: how to share the benefits?

\author[tue]{Loe Schlicher}
\ead{l.p.j.schlicher@tue.nl}
\address[tue]{Eindhoven University of Technology, Department of Industrial Engineering and Innovation Sciences \\ P.O. Box 513, 5600 MB, Eindhoven, The Netherlands }

\author[hec]{Virginie Lurkin\corref{cor1}} 
\ead{virginie.lurkin@unil.ch}
\cortext[cor1]{Corresponding author}
\address[hec]{University of Lausanne, HEC Lausanne Faculty of Business and Economics \\
Quartier de Chamberonne CH-1015, Lausanne, Switserland}

\begin{abstract}
\emph{Horizontal agreements can fall within the scope of exemptions to antitrust competition if they are expected to create pro-consumer benefits. Inspired by such horizontal agreements, we introduce a cooperative game in which a set of transport operators can collectively decide at what price to offer sustainable urban mobility services to a pool of travelers. The travelers choose amongst the mobility services according to a multinomial logit model, and the operators aim at maximizing their joint profit under a constant market share constraint. After showing that various well-known allocation rules (i.e., proportional rules and the Shapley value) do not always generate core allocations, we present a core-guaranteeing allocation rule, the market share exchange rule. This rule first allocates to each transport operator the profit he or she generates under collaboration, and then subsequently compensates those transport  operators  that  lose  part  of  their  market  share,  which  is  paid  by  the  ones  that receive some extra market share. This exchange of market share is facilitated by a unique price,  which can be expressed as the additional return by cooperating per unit of market share. Finally, we show that, under  some  natural conditions, the market share exchange rule still sustains the collaboration when the transport operators need to pay back  part  of  the  joint  profit  to  society.}
 \end{abstract}

\begin{keyword}
game theory, choice-based pricing, cooperative game, core, allocation rules
\end{keyword}

\end{frontmatter}

\section{Introduction}
%Micromobility as sustainable mobility alternative but tough industry
More and more innovative transport solutions, such as e-bikes and e-scooters, are popping up in our urban streets. Together with other shared, small, and light emission-free vehicles, these solutions are known collectively as \textit{micromobility}. Micromobility can help cities and service providers to address unsolved transportation challenges related to urban congestion and pollution (\cite{abduljabbar2021role}). Though micromobility have increased in popularity in major cities, many startups still suffer from profitability issues (\cite{fearnley2020micromobility}). Collaboration \textcolor{black}{between these startups} has been acknowledged as an important factor to unlock the sustainable benefits of micromobility (\cite{moller2020micromobility}). 

%\textcolor{red}{[of course, they also collaborate with regulators, but at the end its about cooperation between the operators. Maybe good to emphasize that most. IF OK remove this text.]}

%Article 101: basic rule = no cooperation but exemptions
While various forms of collaboration (e.g., strategic alliances or public private partnerships) are, \textit{in theory}, possible, creating functioning collaborations remains hard \textit{in practice}. Indeed, following Article 101 of the Treaty on the Functioning of the European Union (TFEU), horizontal agreements between competitors (e.g., micromobility startups) are usually perceived as anti-competitive and are heavily fined under the national laws of the EU Members States. Article 101 however permits some exceptions. For instance, horizontal agreements may qualify for exemption if they create sufficient pro-consumer benefits that outweigh the anti-competitive effects (Article 101(3), TFEU). Besides, they should not eliminate the competition in the relevant market, implying that participants should have a small market share and their combined market share should not exceed a specified limit (e.g., $20\%$ in the Netherlands). A recent example of such an exemption stems from the Webtaxi case in Luxembourg, where the competition authority allowed various taxi companies to use a pricing algorithm to determine their taxi prices. While it was acknowledged that the joint use of such an algorithm constituted a situation of price fixing (i.e, a price agreement between competitors), it was decided that the agreement could be exempted since they expected huge pro-consumer benefits, mainly less waiting time, and lower prices, as well as more rides for drivers (\cite{bostoen2018competition}). At the same time, the combined market shares of the taxi operators would remain far below the threshold set by the Luxembourg authorities. 

%Similar to the horizontal agreement in the taxi sector, micromobility startups in major cities such as Paris, London or Amsterdam could also decide to bundle their forces and set prices collectively. Such an initiative could qualify for an antitrust exemption, since it helps mobility startups with their profitability issues, which is so crucial for the continuation and embedding of environmental friendly mobility solutions in our society. However, according to Article 101 (TFEU), such a horizontal agreement is only allowed if the combined market share does not exceed a specified limit. If the number of participating startups is large enough, this boils down to re-balancing the combined market shares over the startups. Although such a setting will increase the joint profit, it may also lower the market share, and associated profit, of some startups. If these startups are not compensated, it is unlikely that the horizontal agreement will continue to exist. Indeed, for this type of horizontal collaboration, the development of a financial (compensation) scheme, which identifies how to allocate the joint profit amongst the participants, is vital. 

%While these exemptions remain rare in practice (\cite{bostoen2018competition}), they nevertheless create unique potential collaboration opportunities for the micromobility startups that struggle with their profitability issues. More importantly, such exemptions could then help in embedding environmental friendly mobility solutions in our society. 

Similar to the horizontal agreement in the taxi sector, micromobility startups in major cities such as Paris, London or Amsterdam could also decide to bundle their forces and set prices of their mobility services collectively. Such an initiative could also qualify for an antitrust exemption, since it helps mobility startups with their profitability issues, which is so crucial for the continuation and embedding of environmental friendly mobility solutions in our society. However, according to Article 101 (TFEU), such an agreement is only allowed if the combined market share does not exceed a specified limit. To fall within this requirement, the mobility startups could decide to set new prices, while keeping the combined market share stable. Basically, this boils down to a setting where some startups are losing market share, some others are winning market where, but where the overall joint profit (i.e., the sum of the profits of all startups) is increasing. To sustain such form of collaboration, it is important that the losing startups are financially compensated by the winning startups. This can be facilitated by (the development of) an allocation rule, which identifies how to properly distribute the overall joint profit amongst the mobility startups.

Inspired by this joint profit (re)allocation problem, we study a setting in which a set of transport operators (e.g., micromobility startups) can collaborate and decide at what price to offer sustainable urban mobility solutions (e.g., electric scooters or bikes) to a pool of travelers. To better reflect the decisions of these travelers, we assume that they choose among the services offered according to a multinomial logit model, one of the most widely-used disaggregate demand model (\cite{ben2003discrete}). By considering a choice model on the demand side, our setting involves pricing decisions that better capture the supply-demand interactions between the operators objective of maximizing their expected revenue and the travelers objective of maximizing the expected utility (\cite{sumida2019revenue}). To be in line with the conditions associated with the horizontal agreement exemptions, we assume that the transport operators set their prices in such a way that the total joint profit is maximized and their total market share remains constant, and as such remains below the authorized limit. We then formulate a cooperative game for this setting, the transport choice (TC) game, and introduce various intuitive allocation rules. We study to which extent these allocation rules produce core allocations, i.e., allocations that sustain the collaboration because they give no reason for any transport operator, nor group of operators, to break from the collaboration. We show that intuitive proportional rules as well as the well-known Shapley value produce allocations that do not always belong to the core (i.e., these allocations can not sustain the collaboration). We also study a market share exchange allocation rule. This rule first allocates to each transport operator the profit it generates under collaboration. Subsequently, the rule compensates those transport operators that lose part of their market share, which is paid by the ones that receive some extra market share. This exchange of market share is facilitated by a unique price, which can be expressed as the additional return by cooperating per unit of market share. We prove that the allocations of this allocation rule always belong to the core. Finally, we study a setting where the transport operators need to pay back part of the joint profit to society. We show that, under some natural conditions, the allocations of the market share exchange allocation rule are still in the core.

The rest of the paper is organized as follows. Section \ref{LR} provides an overview of the main advancements in the two main disciplines related to this paper: choice-based pricing and cooperative game theory. In Section \ref{preliminaries}, we introduce preliminaries on discrete choice theory and cooperative game theory. Transport choice situations are introduced in Section \ref{sec:transportchoicesituations}. The associated transport choice game is introduced in Section \ref{TCgame}. We study allocation rules for our game in Section \ref{allocationrules}. Then in Section \ref{extension} we study an extended TC game, where players need to pay back part of the joint profit to society. We conclude this paper with final remarks in Section \ref{conclu}. While complete proofs of lemmas and theorems are relegated to Appendix A, a sketch of proof is given in the main text for the theorems that constitute our main results.

\section{Related literature}
\label{LR}

%Pricing problems are important and popular, also in transportation
The determination of optimal prices for different services (or products) is an essential component of operations management \citep{sun2021value}. This is also a delicate task in many organization, as it affects corporate profitability and market competitiveness. The higher the service price, the better the company can cover its costs and generate a profit, but the higher the service price the least attractive the service becomes for the customers \citep{tawfik2018pricing}. This is especially true when newly competitive markets are emerging, which has happened in the transport sector in recent years. As a result of their practical significance, pricing problems have attracted a lot of attention in many fields, including transportation (\cite{azadian2018service,arbib2020competitive,zhong2021lexicographic}). 

%Pricing problems use demand data as input
Most (if not all) pricing problems require demand data as input. In many pricing problems, aggregate representations of demand is used. This aggregate modeling approach is not able to capture the causal relationship between the pricing decisions and the individual customer purchase decisions. To better represent the supply-demand interactions, dis-aggregate demand models have been integrating within pricing problems. The state-of-the-art for the modeling of dis-aggregate demand relies on Discrete-choice modeling (DCM) \citep{bierlaire2017introduction}. 
%Transport researchers have used DCM for more than 40 years, from the pioneer work of \cite{mcfadden1974frontiers} to more recent studies on willingness to pay for self driving vehicles \citep{daziano2017consumers} or willingness to travel with green modes in the context of shared mobility \citep{li2019integrated}. 
Pricing models are usually referred to as \textit{choice-based pricing} if customer's choice behavior is modelled using DCM. These models are  mathematically complex since they are nonlinear and non-convex in prices. Still, as discussed in the next section, the operations research community put remarkable efforts in solving these models because they better reflect the trade-off between the business objective of maximizing the expected revenue and the customer objective of maximizing the expected utility \citep{sumida2019revenue}.

\subsection{Choice-based pricing} 

%choice based pricing
\cite{hanson1996} pioneer this research by showing that the expected revenue function is not concave in prices, even for the simple \textit{multinomial logit} (MNL) model. Subsequent authors have showed that, under uniform price sensitivities across all products, the expected revenue function is concave in the choice probability vector \citep{song2007demand,dong2009dynamic,zhang2013assessing}. \cite{li2011pricing} proved that the concavity remains for asymmetric price-sensitivities, for both the MNL model and the \textit{nested logit} (NL) model that generalizes the MNL model by grouping alternatives into different nests based on their degree of substitution.

%Unique prices LR 
Parallel works have used first-order conditions to show that, under restrictive assumptions on the degree of asymmetry in the price sensitivity parameters, there exist unique price solutions for some logit models. This was demonstrated for the MNL model (e.g., \citep{aydin2000product,hopp2005product,maddah2007joint,aydin2008joint,akccay2010joint}), the NL model (e.g., \cite{aydin2000product,hopp2005product,maddah2007joint,aydin2008joint,akccay2010joint,gallego2014multiproduct,huh2015pricing}), and the \textit{paired combinatorial logit} (PCL) model \citep{li2017optimal}. Lately \cite{zhang2018multiproduct} showed that this result actually holds for the entire family of \textit{generalized extreme value} (GEV) models. This stream of research also includes studies in which pricing decisions are optimized jointly with other decisions such as assortment or scheduling decisions (e.g., \cite{du2016optimal,jalali2019quality,bertsimas2020joint}). 

%ML LR 
Both \cite{gilbert2014mixed} and \cite{li2019product} consider a pricing problem under a \textit{mixed logit model} (ML), a popular choice model that allows the price sensitivity parameter to vary across individuals. Under ML assumption, the concavity property with respect to the choice probabilities breaks down (even for entirely symmetric price sensitivities). The theoretical results obtained for the other logit models therefore do not apply to ML-based pricing problems.  \cite{li2019product} assume a finite number of market segments, with product demand in each segment governed by the MNL model. To solve this problem, the authors propose an algorithm that converges to a local optimum by solving two concave maximization problems, which work as lower and upper bounds for the objective value of the revenue function.  \cite{gilbert2014mixed} consider a ML demand model within a revenue-maximizing network pricing problem. Unlike \cite{li2019product}, the price sensitivity parameter is distributed across the population according to a continuous random variable. To solve this complex problem, the authors rely on a tractable approximation of the ML-pricing problem.  

%Moving from Monopoly to Oligopolies and thus Game Theory
Apart from \cite{li2011pricing}, all above studies consider a monopoly setting. However, in practice multiple groups of decision-makers are simultaneously involved within transport markets. As such, game theory is a suitable framework to analyze these choice-based problems \citep{adler2020review}. Choice-based pricing problems have therefore also been studied from a non-cooperative game theory perspective, as highlighted in the next section.

\subsection{Non-cooperative game theory}

%Competitive setting

In non-cooperative games based on choice-based pricing, important research efforts are made on showing conditions for existence and uniqueness of Nash equilibria. Existence conditions for Nash equilibria for non-cooperative games under MNL and NL models are provided by \cite{milgrom1990rationalizability}, \cite{bernstein2004dynamic} and \cite{li2011pricing}, among others. \cite{aksoy2013price} identify conditions on price bounds and segment market shares that guarantee the existence and uniqueness of equilibrium for a logit-based game involving differentiated firms offering a single-product at a unique price to groups of homogeneous customers. 
Following several empirical evidence, and using the adjusted markup as a single decision variable, \cite{gallego2014multiproduct} show that under mild conditions a unique Nash equilibrium exists for a market with homogeneous demand and nested logit models.

\cite{lin2009dynamic} show the existence of an equilibrium for a dynamic logit-based price competition game between firms selling a single product in a market of substitutable products. The authors propose policies to find the equilibrium in case of full and partial information. In a similar spirit, \cite{levin2009dynamic} consider a stochastic dynamic game based on a generalized choice model of demand where customers are subdivided into market segments. Under certain assumptions with respect to information and competition, the authors show the existence of subgame equilibrium solution for each period of this game.

\cite{morrow2011fixed} present necessary stationarity conditions and analyze numerical methods to compute equilibrium prices for a market with multi-product offer and homogeneous prices under a general ML model of demand. The authors then acknowledge that determining existence or uniqueness of equilibrium prices with general discrete choice models, heterogeneous multi-product firms and heterogeneous consumers is an open problem. In \cite{bortolomiol2021simulation}, a simulation-based heuristic framework is presented to solve pricing problems based on advanced logit models such as the ML, with heterogeneous population, multi-product offer by suppliers and price differentiation. The flexibility of the methodology regarding the choice of the demand model however comes at the expense of pure equilibrium conditions.

In our paper, we also consider a choice-based pricing problem involving multiple firms. Our study therefore fits in the established literature on price optimization under logit choices.
However, inspired by horizontal agreement exemptions (as discussed in the introduction), we assume that the firms can collaborate and collectively decide at what price to offer their services. Cooperative game theory is then the most appropriate methodology to adopt to allocate the associated joint profit between the firms, as explained in the next section.

%fair allocations of the  The next section shows how cooperative game theory has been successfully applied to various type of real-life collaborative settings. 

%The next section shows how cooperative game theory has been successfully applied to various type of real-life collaborative settings. 

%However, inspired by the development of collaborative games over the last decade, our study departs from the competition context and assumes that the firms can collaborate and collectively decide at what price to offer their services. The next section shows how cooperative game theory has been successfully applied to various type of real-life collaborative settings. 

%(\citet{curiel2013cooperative})

%By making use of cooperative game theory, we investigate the stability of such form of collaboration.

\subsection{Cooperative game theory}

Cooperative game theory was successfully applied to various types of real-life collaborative settings. 
For instance, it has been used to identify fair prices for vaccine exchange between countries in times of pandemics (\citet{westerink2020core}), to help museums to decide how to share the profits arising from a museum pass (\citet{ginsburgh2003museum}), to identify fair prices to share railway equipment amongst railway contractors (\citet{schlicher2017probabilistic, schlicher2018pooling,schlicher2020core}), and to help service operations in factories to divide cost savings when they decide to optimally re-balance their production lines. (\citet{anily2017line}). 

The transportation industry also offers many situations suitable for the application of game theory. One can, for example, think of various transport operators that decide to team up to perform (parts of) their logistics operations jointly. By exchanging transport requests among each other, logistics operations can potentially take place in a more efficient and sustainable way. Examples of such application can be found in \citet{lozano2013cooperative}, \citet{engevall2004heterogeneous}, \citet{hezarkhani2016competitive} and \citet{kimms2016core}.

Our paper is enrolled in the line of this existing literature on cooperative games inspired by real-life settings in transport. In particular, it tackles a sensitive topic: price collaboration. To the best of our knowledge, we are the first in this stream of research to focus on collaborative price setting between transport operators. This is not surprising since price fixing is an obvious horizontal agreement and is therefore perceived as anti-competitive and heavily fined under the national laws of the EU Members States. As such, the industry has very little or no motivation to study profit allocation aspects under collaborative price setting. However, as pointed out in our introduction, under some strict conditions, antitrust rules can be repealed, allowing competitors to collaborate on prices. These exemptions have been our source of inspiration for the development of a cooperative game in which a set of transport operators can collectively decide at what price to offer sustainable urban mobility services to a pool of travelers whose purchase decisions are characterized by the MNL model. 

In line with the recent literature on cooperative games inspired by real-life settings in transport, we also investigate the non-emptiness of the core of our game and study fairness properties of various allocation rules (e.g., proportional rules or the Shapley value).

%Examples of cooperative game theory stem from various domains. Among others, they have been successful used

%In the service logistics industry, \citet{schlicher2018pooling,schlicher2020core} showed how cooperative game theory can be helpful to allocate cost savings amongst railway contractors that decide to pool critical assets (e.g., tamping machines) and in a manufacturing setting, \citet{anily2017line} illustrated how cooperative game theory can assist in allocating cost savings amongst service operators that need to (re)balance their production lines.

%For a review on cooperatives games inspired by real-life settings in transport, we refer to \citet{fiestras2011cooperative} and \citet{guajardo2016review}.

%\citet{lozano2013cooperative} used the theory of cooperative games to allocate cost savings that different transport companies (such as shippers)  achieve when they merge their transportation requirements. \citet{engevall2004heterogeneous} illustrated how to allocate travelling costs when several customers are served by one (or several) truck deliveries. In a urban transport setting, \citet{hezarkhani2016competitive} gave a cooperative game solution for gain sharing in consortia of logistic providers where joint planning of truckload deliveries enables the reduction of empty kilometers.

There are some evident similarities between our TC game and the recent contributions of \cite{lardon2019coalitional,lardon2020core} that study cooperative Bertrand oligopoly games. Like us the authors are interested in a set of firms that need to set prices, and have an associated demand function that describes how many customers will opt for that firm. Under the assumption that their demand function is linear in price, they show that their game has a non-empty core, meaning that there are incentives for the firms to cooperate on prices. Our TC game significantly departs from these works by using a non-linear demand function: the MNL model that represents the customers' purchase decisions. Besides, we assume that the sum of the demand functions is stable in price, while it is not in  \cite{lardon2019coalitional,lardon2020core}.

%In Section \ref{TCgame}, we detail this game and in Section \ref{allocationrules} we introduce a core allocation rule for it.  

Finally, it is worth nothing that in a more theoretical context, our TC game also has some similarities with market games (\cite{shapley1969market}. In these games, each player is associated with a set of resources and a concave profit function, identifying the amount of profit realized for a given set of resources. Players collaborate by reallocating their resources to maximize the sum of the utility functions. In our TC game, players are also reallocating resources, namely market share, and are equipped with an implicitly defined utility function that describes to profit per player, for a given amount of market share. However, opposed to market games, in our game the market share is strictly positive for any combination of prices. This due to the logit-form of the MNL model. As such, it is not possible to translate our TC game into a market game, directly. We would like to mention that, although this translation cannot be made directly, we don't rule out that it is still possible. However this transformation may be a difficult as studying the game itself (\citet{anily2017line})

\section{Preliminaries on Discrete Choice and Cooperative Game Theory}
\label{preliminaries}

\subsection{Discrete Choice Theory} \label{subsec:DCT}

Rooted in microeconomics, DCM are powerful operational tools that aim at capturing the causality between a set of explanatory variables and the behavioral choice of economic actors.  The set of alternatives considered as a potential choice is assumed to be finite and discrete, and is referred to as the choice set, denoted $C$ (with $C \not = \emptyset$). For example, in the context of the choice of a transport mode to commute to work, the choice set of an individual could include the car, train, bus, walking, and biking options (i.e., $C = \{\mbox{car, train, bus, walking, biking}\})$.

A fundamental assumption behind these models is that each individual is a rational utility maximizer. It means that, when making a choice among a set of available alternatives, the individual, or the decision maker, is choosing the alternative that maximizes a utility function. The exact specification of this utility function is unknown and therefore typically modeled as a continuous random variable. In this paper, we assume a homogeneous population and that there is no discrimination among individuals. The utility associated with a specific alternative $i \in C$ is thus the same for each individual and given by
\begin{equation*}\label{Eq:utility}
U_{i} = V_{i}(x_{i};\beta) + \varepsilon_{i},
\end{equation*}

\noindent where $V_{i}(x_{i};\beta)$ is a deterministic function of the attributes of the alternatives $x_{i} \in \mathbb{R}$ (e.g., the price of the transport mode), $\beta \in \mathbb{R}$ is a vector of estimable parameters for alternative $i$ (e.g., the willingness to pay), and $\varepsilon_{in}$ is a continuous error term, capturing the specification and measurement errors. The choice model that predicts the probability for an individual to choose alternative $i \in C$ is therefore probabilistic and defined as

\begin{equation*}\label{Eq:choice model}
\mathbb{P}(i) = \mathbb{P}(U_{i} \geq U_{j}, \forall j \in C).
\end{equation*}

In the context of discrete choice, it is custom to assume that all error terms are independent, identically, and extreme value distributed with location parameter $0$ and scale parameter $1$ (i.e., $\varepsilon_{i} \sim$ EV(0,1)). These assumptions lead to a widely used choice model in practice, the logit model whose choice probability is given by

\begin{equation*}\label{Eq:logitprobability}
\mathbb{P}(i) = \frac{e^{ V_{i}(x_{i};\beta)}}{\sum_{j \in C} e^{ V_{j}(x_{j};\beta)}} \mbox{ for all } i \in C.
\end{equation*}

DCM are commonly used in the scientific literature and in practice to understand and predict individual human behavior. The problem then consists in finding the parameters values, i.e., the coefficients of the variables in the utility functions, that maximize the probability that the choice model correctly predicts all observed choices (called the likelihood). In the last decade these choice models have also been more and more used within optimization models to represent more realistically the complexity of human behavior.

\subsection{Cooperative Game Theory} \label{subsec:CGT}

Cooperative game theory primarily deals with the modelling and analysis of situations in which a group of players can benefit from coordinating their actions. In this paper, we model and analyze a specific type of cooperative game: a cooperative game with transferable utility. In what follows, we formally introduce this type of game and discuss desirable properties they may satisfy. We conclude with a description of an allocation rule. 

 A cooperative game with transferable utility, shortly called a (TU) game, is a pair $(N,v)$ where $N$ is a non-empty, finite player set and $v : 2^N \to \mathbb{R}$ a characteristic function with $v(\emptyset) = 0$. 
 A subset $M \subseteq N$ is a coalition and $v(M)$ is the worth coalition $M$ can achieve by itself. This worth can be transferred freely amongst the players. The set $N$ is called the grand coalition. A game $(N, v)$ is monotonic if the value of every coalition is at least the value of any of its subcoalitions, i.e., $v(M) \leq v(K)$ for all $M, K \subseteq N$  with  $M \subseteq K$. When the value of the union of any two disjoint coalitions is larger than or equal to the sum of the values of these disjoint coalitions, a game $(N,v)$ is superadditive, i.e., $v(M) + v(K) \leq v(M \cup  M)$ for all $M, K \subseteq N$ with $M \cap  K = \emptyset$. A game $(N,v)$ is convex if the marginal contribution of any players to any coalition is less than his marginal contribution to a larger coalition, i.e., $v(K \cup \{i\} - V(K)) \geq v(M \cup \{i\}) - v(M))$ for all $M \subseteq K \subseteq N \backslash \{i\})$ and all $i \in N$.

An allocation for a game $(N, v)$ is an $N$-dimensional vector $x \in \mathbb{R}^N$ describing the payoffs to the players, where player $i \in N$ receives $x_i$. An allocation is called efficient if $ \sum_{i \in N} x_i = v(N)$. This implies that all worth is divided amongst the players of the grand coalition $N$. An allocation is individual rational if $x_i \geq v(\{i\})$ for all $i \in N$ and stable if no group of players has an incentive to leave the grand coalition $N$, i.e. $ \sum_{i \in M} x_i \geq v(M)$ for all $M \subseteq N$. The set of efficient and stable allocations, called the core of $(N, v)$, is denoted by
$$\mathscr{C}(N,v):= \left\{ x \in \mathbb{R}^N  \hspace{2mm} \bigg \vert \hspace{2mm} \sum_{i \in M} x_i \geq v(M) \mbox{ for all } M \subseteq N \mbox{ and } \sum_{i \in N} x_i = v(N)\right\}.$$

An allocation rule is a function $f$ that assigns to any game $(N,v)$ in a class of cooperative games a vector $f(N,v) \in \mathbb{R}^N$ satisfying $\sum_{i \in N} f_i(N,v) = v(N)$. A well-known allocation rule defined on the set of all games is the Shapley value (\citet{shapley1953value}). This allocation rule assigns to each player a weighed average over all marginal contributions (s)he can make to any possibly coalition. Formally, for any game $(N,v)$ the Shapley value can be defined by:
$$ SV_i = \sum_{M \subseteq N \backslash \{i\}} \frac{ \vert M \vert ! (\vert N \vert - 1 - \vert M \vert )!}{\vert N \vert !} (v(M \cup \{i\}) - v(M)) \mbox{ for all }i \in N.$$

\section{Transport Choice situations} \label{sec:transportchoicesituations}

We consider a setting in which a group of homogeneous travelers is buying mobility services from a set of $N \subseteq \mathbb{N}$ transport operators. Each operator $i \in N$ offers one micromobility service (e.g., a e-bike or a segway) against price $p_i \in \mathbb{R}_+$ and cost price $c_i \in \mathbb{R}_+$. The mobility choices of travelers are represented using the logit model described in Section \ref{subsec:DCT}. In doing so, we let, per mobility service $i \in N$, function $V_i$ be defined as:
%In particular, for each mobility service offered by a transport company, we assume a linear-in-parameters utility specification that includes an alternative-specific constant, $\alpha_i$, and the price of the mobility service $p_i$, weighted by the price sensitivity parameter $\beta$:
\begin{equation*}\label{Eq:deterministicutility}
V_{i} = \alpha_i -\beta p_i,
\end{equation*} 

where $\alpha_i \in \mathbb{R}_+$ is an alternative-specific constant and $\beta \in \mathbb{R}_+$ a price sensitivity parameter. Moreover, the choice set $C$ of the logit model consists of ($i$) the transport operators, where each operator refers to the operator-specific mobility service, and ($ii$) the choice to not buy any service, for which the deterministic utility is normalized to zero. We denote the choice of not buying service by $0$, and consequently $C= N \cup \{0\}$ and $V_0 =0$.

By using the choice probabilities of the logit model, we can define the share of \textcolor{black}{travelers} that opts for mobility service $i \in N$ as follows
\begin{align}
\label{eq:MNLproba} 
 \frac{e^{\alpha_i - \beta p_i}}{\sum_{j \in C} e^{\alpha_j - \beta p_j}} = \frac{e^{\alpha_i - \beta p_i}}{1 + \sum_{j \in N} e^{\alpha_j - \beta p_j}}.\end{align}
\textcolor{black}{From now on, we refer to this share of travelers as the market share of transport operator $i$.} Given this market share, the profit of transport operator $i \in N$ is defined by:
\begin{align}
\label{eq:MNLprofit} (p_i - c_i) \cdot \frac{e^{\alpha_i - \beta p_i}}{1 + \sum_{j \in N} e^{\alpha_j - \beta p_j}}.
\end{align}
We summarize this setting by tuple $\theta = (N,p,c,\alpha,\beta)$ with $N$ the set of transport operators, $p = (p_i)_{i \in N}$ the vector of prices, $c=(c_i)_{i \in N}$ the vector of cost prices,  $\alpha = (\alpha_i)_{i \in N}$ the vector of alternative-specific constants, and $\beta$ the price sensitivity parameter. We refer to $\theta$ as the transport choice situation and let $\Theta$ be the set of all possible transport choice situations. \\

We now illustrate our TC situation with a fictitious example. 

\begin{example} \label{firstexample1} Let $\theta \in \Theta$ with $N=\{1,2,3\}$, $p = (6,8,15)$, $c = (8,4,1)$, $\alpha = (1,0.5,1.5)$ and $\beta= 0.36$. The prices, corresponding market shares and associated profits of the transport operators are presented in Table \ref{monopolyprices}. $\hfill \diamond$

%Note that, while this example is completely fictitious, operator 1 could represent a loss-making startup whose strategy is to first drive towards capturing market share.

\begin{table}[h!] \centering
\begin{tabular}{cccc} \hline
    $i$ & 1 & 2 & 3 \\ \hline
    Price $i$ & 6.0 & 8.0 & 15.0 \\
    market share $i$ & 0.220 & 0.065 & 0.014 \\
  profit $i$ & -0.440 & 0.260 & 0.199 \\  \hline
\end{tabular}
\caption{Prices, market shares and profits of the transport operators of situation $\theta$}
\label{monopolyprices}
\end{table} 
\end{example}

%0.1861233515793806, 0.050724530731428526, 0.008384708535466835
%-0.3722467031587612, 0.2028981229257141, 0.11738591949653569

\begin{remark}
If the aim of the transport operators is to maximize their profits, it is reasonable to consider only those $\theta \in \Theta$ for which $p$ is a Nash equilibrium (i.e., vector $p$ is such that no transport operator $i \in N$ would unilaterally deviate his or her $p_i$, for the given $N$, $c$, $\alpha$ and $\beta$). It can be shown (see Appendix A), that for those situations, $p$ satisfies
\begin{equation}\label{NEprices}
p_i = \frac{1+W\left(\frac{e^{\alpha_i-1-\beta c_i}}{1+ \sum_{j \in N \backslash i} e^{\alpha_j - p_j\beta} }\right)}{\beta} + c_i \hspace{2mm} \mbox{ for all } i \in N,
\end{equation} 

\noindent where $W$ denotes the Lambert $W$ function.

%Please note that all subsequent theoretical results in this paper applies for any $\theta \in \Theta$, therefore also including the $\theta \in \Theta$ for which $p$ is a Nash equilibrium.

\end{remark}

Rather than operating individually, the transport operators can decide to cooperate. If the transport operators decide to do so, they will set new prices, which maximize the sum of the operator-specific profits. By doing so, they have to take into account that the total market share remains stable\footnote{\textcolor{black}{With this assumption, we guarantee that the operators cannot start dominating the market, which is one of the governmental rules that needs to be satisfied in order to allow for an antitrust law exception.}}. Formally, if the transport operators in $N$ decide to collaborate, they are \textcolor{black}{facing the following non-linear\textcolor{black}{, constrained,} optimization problem:}

\begin{equation} \label{optimizationproblem}  \begin{aligned} \mathscr{P}:= & \hspace{5mm}  \max_{x \in \mathbb{R}^N} \sum_{i \in N } (x_i - c_i) \frac{e^{\alpha_i - \beta x_i}}{1+\sum_{j \in N} e^{\alpha_i - \beta x_j}} \\& \mbox{ s.t.} \sum_{i \in N} \left( \frac{e^{\alpha_i - \beta p_i}}{1+\sum_{j \in N}e^{\alpha_j - \beta p_j}} \right) = \sum_{i \in N} \left( \frac{e^{\alpha_i - \beta x_i}}{1+\sum_{j \in N}e^{\alpha_j - \beta x_j}} \right)\end{aligned} \end{equation}

%\begin{equation} \label{optimizationproblem}  \begin{aligned}  & \hspace{10mm} \mathscr{P} = \max_{x \in \mathbb{R}^N} \sum_{i \in N } (x_i-c_i) \frac{e^{\alpha_i - \beta x}}{1+\sum_{j \in N} e^{\alpha_j - \beta x }} \\& \mbox{ s.t.} \sum_{i \in N} S_i = \sum_{i \in N} \frac{e^{\alpha_i - \beta x_i}}{1 + \sum_{j \in N} e^{\alpha_j - \beta x_j}}.\end{aligned} \end{equation}

%where $S_i(x) = \frac{e^{\alpha_i-\beta x_i}}{\sum_{j \in N} e^{\alpha_j - \beta x_j}}$ and $\pi_i(x) = (x_i-c_i) \cdot S_i(x)$ for all $x \in \mathbb{R}^N_{\geq 0}$ and all $i \in N$. 

\textcolor{black}{We refer to an optimal solution of $\mathscr{P}$ as an optimal price vector $p^*$ and call the optimal value of $\mathscr{P}$ the optimal joint profit $\mathscr{P}^*$.} There exists a closed-form expression for $p^*$ and $\mathscr{P}^*$. Before presenting them, we first introduce some new notation. For each TC situation $\theta \in \Theta$, we define $D(x) = \sum_{i \in N} e^{\alpha_i - \beta x_i}$ for all $x \in \mathbb{R}^N$. Note that this new notation can be used to describe the total market share compactly (i.e., $\frac{\sum_{i \in N}e^{\alpha_i - \beta p_i}}{\sum_{j \in N} e^{\alpha_i - \beta p_j} + 1} = \frac{D(p)}{D(p) + 1}$). Now, we are ready to present an optimal price vector $p^*$ and the associated optimal joint profit $\mathscr{P}^*$.

\begin{theorem}\label{optimres} For each TC situation $\theta \in \Theta$ an optimal price vector $p^*$ is given by
$$ p_i^* = c_i + \frac{1}{\beta} \ln \left( \frac{D(c)}{D(p)} \right) \mbox{ for all } i \in N,$$
and the associated optimal joint profit equals $\mathscr{P}^* = \frac{D(p)}{\beta(D(p)+1)} \ln \left( \frac{D(c)}{D(p)}\right).$
\end{theorem}

The proof of Theorem \ref{optimres} consists of three steps. First we identify how our optimization problem $\mathscr{P}$ relates to another optimization problem. This optimization problem has a much simpler form of constraint. Then, we identify an optimal price vector and the associated optimal value for this optimization problem, by using a Lagrangian type of optimality result from \cite{bazaraa2013nonlinear}. Finally, we relate back these outcomes to $\mathscr{P}$. \bigskip

We now make some remarks regarding Theorem \ref{optimres}.

\begin{remark} \textcolor{black}{Observe that vector $p^*$ is only player-specific in the cost price. So, if $c_i=0$ for all $i \in N$, every transport operator will select the same price, under collaboration.} \end{remark}

\begin{remark} \label{remark2} The gain of collaboration is always non-negative, i.e., 
$$\mathscr{P}^* - \sum_{i \in N} (p_i - c_i) \cdot \frac{e^{\alpha_i - \beta p_i}}{D(p)+1} \geq0,$$ 
which is due to the fact that price vector $p$ is a feasible solution of $\mathscr{P}$. 
\end{remark}

\begin{remark} \label{remark4} If $p_i -  c_i =p_j -c_j$ for all  $i,j \in N$, the gain of collaboration equals zero, i.e.,
$$ \mathscr{P}^* - \sum_{i \in N} (p_i - c_i) \cdot \frac{e^{\alpha_i - \beta p_i}}{D(p) + 1}= \frac{(p_j-c_j) D(p)}{D(p)+1} -  \frac{(p_j-c_j) D(p)}{D(p)+1} = 0 \mbox{ for all } j 
\in N.$$
   \textcolor{black}{Note, this condition describes a TC situation with a constant marginal profit per operator and a total market share of $D(p)/(D(p) + 1)$. As such, it can be used as a benchmark for situations with the same total market share, but where operators gain from collaboration.} 
    \end{remark}
    
%    We will refer to such TC situations as gain-neutral TC situations and denote them by $\theta^{\Delt}$.

%we want to compare our situation with the one in which the total market share is the same and where collaboration does not lead to additional gains. From Remark X, we learned that such situations occur when, for instance, if $p_i = c_i + 1/\beta$ for all $i \in N$. We also learned that in such a setting, the joint profit equals $\mathscr{P}^* = \frac{D^N(p)}{\beta(D^N(p)+1)}$.   

We now illustrate how Theorem \ref{optimres} applies to our TC situation of Example \ref{firstexample1}.

\begin{example} \label{secondexample} Reconsider the situation of Example \ref{firstexample1}. In Table \ref{nashprices}, we present an optimal price vector and the corresponding market share\footnote{Note that because of rounding the total market share (0.008 + 0.024 + 0.214) seems to have changed, but this is not the case. The total market share still remains 0.245.} and profit per transport operator. $\hfill \diamond$ 
\begin{table}[h!] \centering
\begin{tabular}{cccc} \hline
    $i$ & 1 & 2 & 3 \\ \hline
    optimal price $i$ & 13.980 & 9.980 & 6.980 \\
    market share $i$ & 0.012 & 0.032 & 0.255 \\
    profit $i$ & 0.074 & 0.190 & 1.523 \\  \hline
\end{tabular}
\caption{Prices, market share and profits of the transport operators of situation $\theta$}
\label{nashprices}
\end{table} 
\end{example}

%(array([13.90487637,  9.90487637,  6.90487637]),
% [0.007881033778463145, 0.023675933910833172, 0.21367562315697952],
 %[0.04653653014961818, 0.1398034627478418, 1.2617281385669077])

From Table \ref{monopolyprices} of Example \ref{firstexample1} and Table \ref{nashprices} of Example \ref{secondexample}, we learn that the joint profit, which is $1.787$, exceeds the sum of individual profits without collaboration, namely $-0.440+0.260+0.199 =0.019$. However, at the same time, we also observe that the individual profit of transport operator 2 decreases (from 0.260 to 0.190). So, in case of collaboration among the three transport operators, it would be natural that operator 1 and operator 3 would compensate operator 2 in some way. But, by how much? In the upcoming section, we address this question by making use of cooperative game theory. 

%-0.05196266073651139
%A central question is then: what should be the value of this compensation?

%With cooperative game theory, one can answer such a question. More specifically, this theory aims to determine how to divide the payoff of a coalition between members such that no one would find it worthwhile to deviate from the existing agreement (\citet{adler2020review}). In the next section, we will introduce and analyze such a cooperative game, which arises from our underlying TC situation. 

\section{A Cooperative Transport Choice Game}
\label{TCgame}

In this section, we introduce a cooperative game, associated to our transport choice situation. Formally, for each TC situation $\theta \in \Theta$, we introduce a cooperative game $(N, v^{\theta})$, where $N$ represents the set of players (i.e., transport operators) and $v^{\theta}$ represents the characteristic value function. In this game, $v^{\theta}(M)$ reflects the joint profit coalition $M \subseteq N \backslash \{\emptyset\}$ can realize. This joint profit is obtained by taking into account that ($i$) the sum of the market shares of the players in $M$ remains stable (i.e., the new vector of prices should be such that the sum of their market shares remains the same) and ($ii$) all players outside coalition $M$ (i.e., players in $N \backslash  M$) keep their initially set prices. This game, which we refer to as a cooperative transport choice game, is formally defined as follows.

%\footnote{This assumption is not very restrictive. Later on, we we will show that most of our results remain valid under other assumptions (e.g., players outside $M$ keep their total market share the same, or players outside $M$ will, in anticipation on the prices of $M$, set prices that minimize the profit of coalition $M$}. 

\begin{definition} For every TC situation $\theta \in \Theta$, the associated cooperative transport choice (TC) game $(N,v^{\theta})$ is defined by
    \begin{equation} \begin{aligned} v^{\theta}(M) = \max_{x \in \mathbb{R}^M} \sum_{i \in M }(x_i - c_i) \frac{e^{\alpha_i - \beta x_i}}{1+ \sum_{j \in M} e^{\alpha_j - \beta x_j} + \sum_{j \in N \backslash M}  e^{\alpha_j - \beta p_j}} \hspace{18mm} \\
 s.t.  \sum_{i \in M} \left( \frac{e^{\alpha_i - \beta p_i}}{1+\sum_{j \in N} e^{\alpha_i - \beta p_j}  }\right) = \sum_{i \in M} \left( \frac{e^{\alpha_i - \beta x_i}}{1+\sum_{j \in M}e^{\alpha_i - \beta x_j} + \sum_{ j \in N \backslash M} e^{\alpha_i - \beta p_j}} \right)\end{aligned} \end{equation} 
 for all $M \subseteq N \backslash \{\emptyset\}$ and $v^{\theta}(\emptyset) = 0$.
\end{definition}
%For this, we assume that ($i$) players in $M$ can select new prices such that the market share of coalition $M$ remains the same and ($ii$) players outside coalition $M$ keep their total market share the same as well. This latter assumption holds for instance true in a setting in which the players outside $M$ for instance keep their old prices. 

%In order to define $v^{\theta}(M)$ formally, we first introduce some new definitions. 

Similar to the optimization problem in Section \ref{sec:transportchoicesituations}, we present a closed-form expression for the optimal joint profit of any coalition $M \subseteq N \backslash \{ \emptyset\}$. Before doing so, we need to introduce some coalition-specific notation. For each TC situation $\theta \in \Theta$ and each $M \subseteq N$, we let $D^M(x) = \sum_{i \in M} e^{\alpha_i-\beta x_i}$. Please, note that we have $D^N(x) = D(x)$ for all $x \in \mathbb{R}$. Now, we are ready to present the closed-form expression for any coalition. 

\begin{theorem} \label{res:lemma} For every TC situation $\theta \in \Theta$ it holds, for all $M \subseteq N \backslash \{\emptyset\}$, that 

$$v^{\theta}(M) = \frac{D^M(p)}{\beta(D^N(p)+1)} \ln \left( \frac{D^M(c)}{D^M(p)}\right)$$

\end{theorem}

The structure of the proof of Theorem \ref{res:lemma} is similar to the structure of the proof of Theorem \ref{optimres}. We now present an example of a TC game.

\begin{example} Reconsider the setting of Example \ref{firstexample1}. The coalitional values of TC game $(N,v^{\theta})$ are represented in Table \ref{game1} below. 
\begin{table}[h!] \centering
\begin{tabular}{ccccccccc} \hline
    $M$ &  $\{\emptyset\}$ & \{1\} & \{2\} & \{3\} & \{1,2\} & \{1,3\} & \{2,3\} & \{1,2,3\}\\ \hline
    $v^{\theta}(M)$ & 0 & -0.440  &  0.260  &  0.199  &  0.230   &  1.485 &  0.756  & 1.787  \\ \hline
    \end{tabular}
\caption{Coalitional values of game $(N,v^{\theta})$}
\label{game1}
\end{table} 

Please, observe that the coalitional values of the individual coalitions (-0.440, 0.260 and 0.199) match with the profits of Table \ref{monopolyprices}, and that the coalitional value of the grand coalition (1.787) matches with the sum of the profits of Table \ref{nashprices}. $\hfill \diamond$
\end{example}

\begin{remark} Some readers may see similarities between our TC game and a market game (see e.g., \cite{anily2017line}). We want to emphasize that it is not straightforward to recognize our TC game as a market game. In a market game, players can freely reallocate their resources, implying that some of the players may end up with no resources at all. However, in our game, where market shares could be recognized as resources, it is not possible to assign a player with no market share (note that $D^{\{i\}}(p)/(D^N(p)+1) >0$ for all $p \in \mathbb{R}^N$). As such, our game does not fall in the framework of a market game, directly. \end{remark}

\begin{remark} \label{delaatsteremarkofzo}
As discussed in the preliminaries (Section 3), it is natural to study a cooperative game on  properties like monotonicity, superadditivity and convexity. As a side result, we would like to share that our TC game is superadditive, but not monotonic nor convex.
\end{remark}

\section{Allocation rules}
\label{allocationrules}

The central question is this section is how players of our TC game should distribute the joint profit, in order to sustain the collaboration. In the game theory literature, it is common to address this question by introducing several allocation rules and by investigating whether their allocations belong to the core (see, e.g., \citet{westerink2020core}). In this paper, we also follow this approach. In particular, we introduce four intuitive allocation rules and study whether their allocations belong the core. Recall (from Section 3.2) that the core is the set of allocations that makes no group of players break from the grand coalition, i.e., for any TC situation $\theta \in \Theta$ and associated $(N,v^{\theta})$ the core is given by

%In literature, core allocations, i.e., allocations that make no group of players break from the grand coalition, are considered as  fair allocations that sustain collaboration (see, e.g. \citet{westerink2020core}. In particular, for every TC situation $\theta \in \Theta$ and associated TC game $(N,v^{\theta})$, the core is given by
$$\mathscr{C}(N,v^{\theta}) = \left\{ x \in \mathbb{R}^N \hspace{2mm} \bigg \vert \hspace{2mm} \sum_{i \in M} x_i \geq v^{\theta}(M) \mbox{ for all } M \subseteq N  \mbox{ and } \sum_{i \in N} x_i = v^{\theta}(N) \right\}. $$ 

For each of the four allocation rules, we generate 10,000 random TC situations $\theta \in \Theta$ for $N=3, 4$, and $5$ players, respectively. For each $\theta \in \Theta$, we generate random numbers for the vector of costs ($c$) as well as for the utility parameters ($\alpha$ and $\beta$). More specifically, for each $i \in N$, cost $c_i$ and parameter $\alpha_i$ are drawn for a uniform distribution in the interval [0.5, 15.0] with a 0.5 step. The price sensitivity parameter $\beta$ is drawn from a uniform distribution in the interval (0,1] with step 0.1. The price vector $p$ is computed based on the cost and utility parameters, as shown in equation (\ref{NEprices}), i.e., vector $p$ is a Nash equilibrium. Per generation of 10,000 TC situations, and per allocation rule, we are then interested in the outcome 
$$\frac{\mbox{\# of TC situations for which allocation is in the core}}{\mbox{total \# of TC situations (10,000)}}.$$

%We also repeat these numerical experiments for TC situations in which the initial price vector contains the Nash equilibrium prices (computed based on cost and utility parameters values, as shown in equation (\ref{NEprices})).

 In the upcoming paragraphs, we introduce our four allocation rules and present and discuss the associated outcomes per allocation rule.

%Let us recall that for all experiments, the optimal prices are given by Theorem \ref{optimres}. 

%Finally, We also repeat the experiments for TC situations in which the initial price vector contains the Nash equilibrium prices (as given in Equation \ref{NEprices}).

%In this section we investigate several, intuitive, allocation rules. In particular, we will generating random TC situations and check whether the allocation rules generate allocations that belong to the core. More specifically, we simulate 10,000 TC situations $\theta \in \Theta$ for $N=$3, 4,and 5 players, respectively. For each TC situation $\theta$, we generate random numbers for the vectors of prices, and cost prices ($p$ and $c$) as well as for the utility parameters ($\alpha$ and $\beta$). More specifically, the prices $p_i$ are drawn from a uniform distribution in the interval [0.5, 20.0] with a 0.5 step. For costs $c_i$ and parameters $\alpha_i$ we use the range [0.5, 5.0], again with a 0.5 step. The price sensitivity parameter $\beta$ is drawn from a uniform distribution in the interval [0,1]. Finally, We also repeat the experiments for TC situations in which the initial price vector contains the Nash equilibrium prices (as given in Equation \ref{NEprices}).

\paragraph{\textbf{Allocation rules 1 and 2: Proportional}} 

%For instance, these rules have shown to be relevant in collaborative settings, such as pooling of spare parts amongst various companies the high-tech industry (see, e.g., \citet{karsten2014pooling, karsten2015resource}) and pooling of tamping machines amongst contractors in the railway sector (see, e.g., \citet{schlicher2018pooling}). 

%In cooperative game theory literature, these rules are well studied, both from a theoretical and practical perspective. For instance, in \citet{karsten2014pooling, karsten2015resource}) and \citet{schlicher2018pooling}, it is shown that such proportional rules admit interesting fairness properties (e.g., their allocations are core guaranteed) and can be applied to settings in which companies can pool resources.

Allocation rules that have a long tradition when costs, profits, or savings have to be shared among different agents, are proportional rules (see, e.g., \citet{moulin1987equal}). As the name suggests, these rules allocate the worth to the players in a proportional way. A common proportional rule is to divide the worth proportional to the value of the individual coalitions. Formally, for any $\theta \in \Theta$ and associated TC game, this individual-proportional allocation rule is given by

%Allocation rules that have a long tradition when costs, profits, or savings have to be shared among different agents, are proportional rules (see, e.g., \citet{moulin1987equal}). As the name suggests, proportional rules allocate the worth of the grand coalition (i.e., $v^{\theta}(N))$, to the different players in a proportional way. A common proportional rule is to divide the value of the grand coalition proportional to the profit of the individual players (i.e., the value of the singleton coalitions). Formally, for any $\theta \in \Theta$ and associated TC game, this individual-proportional allocation rule is given by
$$\mbox{I-PROP}_i = \frac{v^{\theta}(\{i\})}{\sum_{j \in N}v^{\theta}(\{j\})} \cdot v^{\theta}(N) \mbox{ for all } i \in N.$$ 

As an alternative, one can also decide to divide the worth  proportional to the initial market shares of the players. In that case, for every $\theta \in \Theta$ and associated TC game, the allocation rule, which we call the  market share-proportional rule, is given by

%Another proportional rule would be that each player receives a fraction of the grand coalition value proportional to his or her market share. In that case, for every $\theta \in \Theta$ and associated TC game, this market share-proportional allocation rule is given by
$$ \mbox{M-PROP}_i = \frac{S_i}{\sum_{j \in N} S_j} \cdot v^{\theta}(N) \mbox{ for all } i \in N.$$ 

 The results of the experiments with respect to the above proportional rules are presented in Table \ref{proportionalrules}. We can see that almost all M-PROP allocations don't belong to the core. We also see that the I-PROP allocation rule performs much better, with around 95\% of allocations belonging to the core for a 3 player TC game. Still it does not always lead to core allocations. In particular, the number of allocations not belonging to the core increases as the size of the game grows. 
 
% \textcolor{red}{Hence, both rules are not able to properly compensate players for exchanging market shares.}

%\begin{table}[h!] \centering
%\begin{tabular}{lcccccc} 
%      &  \multicolumn{3}{c}{From random prices} &  \multicolumn{3}{c}{From Nash prices}\\ 
%            &  \multicolumn{3}{c}{$N$} &  \multicolumn{3}{c}{$N$}\\ 
%      \cline{2-7}
%     & 3 & 4 & 5 & 3 & 4 & 5  \\ \hline
%    I-PROP & 0.3955 & 0.1674 & 0.0769 & 0.5027 & 0.2856 & 0.1771  \\ 
%    M-PROP & 0.0193 & 0.0043 & 0.0009 & 0.0002 & 0.0000  & 0.0000  \\ 
%    \hline
%\end{tabular}
%\caption{Outcome for the two proportional rules (based on %10,000 TC situations)}
%\label{proportionalrules}
%\end{table} 

\begin{table}[h!] \centering
\begin{tabular}{lcccccc} 
            &  \multicolumn{3}{c}{$N$} \\
      \cline{2-4}
     & 3 & 4 & 5   \\ \hline
    I-PROP & 0.9460  &  0.8959 &  0.8538 \\ 
    M-PROP & 0.0001  &  0.0000 & 0.0000  \\ 
    \hline
\end{tabular}
\caption{Outcome for the two proportional rules (based on 10,000 TC situations)}
\label{proportionalrules}
\end{table} 

We now illustrate these proportional allocation rules to our TC situation of Example \ref{firstexample1} and investigate whether their allocations do belong to the core or not. 

\begin{example} \label{firstexample} Reconsider the TC situation $\theta \in \Theta$ and game $(N,v^{\theta})$ of Example \ref{firstexample1}. %The market shares and associated profits of the players before and after collaborating are presented in Table \ref{beforeafter}. 
The allocations of the proportional allocation rules for $(N,v^{\theta})$ are reported in Table \ref{proportionalalloc}.  $\hfill \diamond$

%\begin{table}[h!] \centering
%\begin{tabular}{lcccccc} 
%      &  \multicolumn{3}{c}{Before collaboration} &  \multicolumn{3}{c}{After collaboration}\\ 
%            &  \multicolumn{3}{c}{$i$} &  \multicolumn{3}{c}{$i$}\\ 
%      \cline{2-7}
%     & 1 & 2 & 3 & 1 & 2 & 3  \\ \hline
%   price & 4.00 & 15.00 & 2.00 & 6.591 & 7.091  & 7.591 \\ 
% market share & 0.224  & 0.203 & 0.451  & 0.173 & 0.447  & 0.258  \\  
%profit & 0.783  & 2.836  & 0.225  & 1.052 & 2.721 & 1.570  \\  \hline
%\end{tabular}
%\caption{Prices, market shares and associated profits before and after collaborating}
%\label{beforeafter}
%\end{table} 

\begin{table}[h!] \centering
\begin{tabular}{lccc} 
& & \multicolumn{2}{c}{$\sum_{i \in M} x_i$} \\
$M$ & $v^{\theta}(M)$  &  I-PROP &  M-PROP \\
\hline
$\{1\}$ & -0.440 & -42.101  & 
1.314 \\
$\{2\}$ & 0.260 &  24.859 & 
0.388 \\
$\{3\}$ & 0.199 &  19.029  & 
0.085 \\
$\{1,2\}$ & 0.230 &  -17.242  & 1.702 \\
$\{1,3\}$ & 1.485 & -23.072 & 
1.399 \\
$\{2,3\}$ & 0.756 &  43.889  & 
0.473  \\
$\{1,2,3\}$ & 1.787 &  1.787   &  
1.787  \\
\hline
\end{tabular}
\caption{Illustration of proportional allocation rules}
\label{proportionalalloc}
\end{table} 

%By allocating the grand coalition worth based on the singleton values, player 3 receives a relatively small proportion of the total profit (less than 6\%), while the new price offered by this player, and the corresponding new market share, results in a significant profit increase that account now for almost 30\% of the grand coalition value.

Recall that an allocation $x \in \mathbb{R}^3$ is in the core if 
\begin{itemize}
    \item $x_1 + x_2 + x_3 = v^{\theta}(N)$
    \item $x_1 \geq v^{\theta}(\{1\})$, $x_2 \geq v^{\theta}(\{2\})$, $x_3 \geq v^{\theta}(\{3\})$,
    \item  $x_1 + x_2 \geq v^{\theta}(\{1,2\})$, $x_1 + x_3 \geq v^{\theta}(\{1,3\})$, and $x_2 + x_3\geq v^{\theta}(\{2,3\})$. 
\end{itemize}

%The stability condition "" is broken. In words, this means that 

I-PROP is not in the core, since $\mbox{I-PROP}_1 + \mbox{I-PROP}_2 = -17.242 < 0.230 = v^{\theta}(\{1,2\})$. That means, players 1 and 2 together can earn more by breaking up and forming a new coalition together. Similarly, M-PROP is not in the core, since $\mbox{M-PROP}_3 = 0.085 <  0.199 = v^{\theta}(\{3\})$. That means, player 3 is better of by working individually.

\end{example}

\paragraph{\textbf{Allocation rule 3: Shapley value}}

Another well-known allocation rule is the Shapley value. This allocation rule is introduced in \citet{shapley1953value}, and has shown to be applicable in various settings, such as cost sharing in horizontal cooperation among shippers (\citet{lozano2013cooperative}), carpool problems (\citet{naor2005fairness}), and data sharing settings (\citet{dehez2013data}). In words, the Shapley value assigns to each player a weighted average over all marginal contributions a player can make to any possible coalition. Formally, for any TC situation $\theta \in \Theta$ and associated TC game, the Shapley value is defined as:
$$ SV_i =  \sum_{M \subseteq N \backslash \{i\}} \frac{ \vert M \vert ! (\vert N \vert - 1 - \vert M \vert )!}{\vert N \vert !} (v^{\theta}(M \cup \{i\}) - v^{\theta}(M)) \mbox{ for all }i \in N.$$

%Next to its applications, the Shapley value is also of interest from a theoretical point of view. For instance, it is the only rule that satisfies the properties efficiency, symmetry, additivity, and null player, simultaneously. However, the Shapley value also has a major drawback: its allocations are not core guaranteed. In particular, if a cooperative game is not convex (see Section 2.2.), the allocations of the Shapley value may fall outside the core. 

Our numerical results in Table \ref{SVrule} demonstrate that also the shapley value does not always belong the core. In terms of performance, the Shapley value is slightly better than the I-PROP rule (see Table \ref{proportionalalloc}). The Shapley value does belong to the core in 96\% of the instances for a 3 player TC game, and in 90\% for a 5 player TC game. 

%Smiliar to the proportional rules, in most cases one of the individual rationality conditions is broken, meaning that at least one player receives less when collaborating than when working alone. 

\begin{table}[h!] \centering
\begin{tabular}{lcccccc} 
            &  \multicolumn{3}{c}{$N$} \\ 
      \cline{2-4}
     & 3 & 4 & 5  \\ \hline
   SV &  0.9620 & 0.9303 & 0.8991 \\  \hline
\end{tabular}
\caption{Outcome for the Shapley value (based on 10,000 TC situations)}
\label{SVrule}
\end{table}

%\begin{table}[h!] \centering
%\begin{tabular}{lcccccc} 
%      &  \multicolumn{3}{c}{From random prices} &  \multicolumn{3}{c}{From Nash prices}\\ 
 %           &  \multicolumn{3}{c}{$N$} &  \multicolumn{3}{c}{$N$}\\ 
  %    \cline{2-7}
  %   & 3 & 4 & 5 & 3 & 4 & 5  \\ \hline
  % SV & 0.4623 & 0.2055 & 0.0984 & 0.6641 & 0.4119 & 0.2665 \\  \hline
%\end{tabular}
%\caption{Outcome for the Shapley value (based on 10,000 TC situations)}
%\label{SVrule}
%\end{table} 

%From  the previous paragraphs, we learned that players should somehow be compensated for exchanging their market shares.

We provide below an example, illustrating the calculation of SV.

\begin{example} Reconsider the TC situation $\theta \in \Theta$ and game $(N,v^{\theta})$ of Example \ref{firstexample1}. Then, $SV = (0.407,	0.392,	0.989)$. One can check that the Shapley does not belong to the core.
\end{example}

\paragraph{\textbf{Allocation rule 4: Market Share Exchange}} From the previous paragraphs, we learned that the allocations of the proportional rules and the Shapley value do not always belong to the core. One reason could be that these allocation rules do not explicitly compensate for the exchange of market share between players. Therefore, in this paragraph, we study an allocation rule that does explicitly compensate for this exchange of market share. In particular, we study an allocation rule that first allocates to each player the profit he/she generates under full collaboration, i.e., player $i \in N$ receives $(p_i^* - c_i) \frac{D^{\{i\}}(p^*)}{D^N(p^*)+1}$. Thereafter, we identify for each player $i \in N$ the increase (or decrease) in the market share, which is $\left( \frac{D^{\{i\}}(p^*)}{D^N(p^*) + 1} - \frac{D^{\{i\}}(p)}{D^N(p) + 1}\right)$. Player $i$ then receives a price $\phi$ for each exchanged unit of market share, and pays the same price for each extra unit of market share. Formally, for each $\theta \in \Theta$ and associated TC game, the market share exchange (MSE) rule is given by
$$ MSE_i = (p_i^* - c_i) \frac{D^{\{i\}}(p^*)}{D^N(p^*)+1} -  \phi \left( \frac{D^{\{i\}}(p^*)}{D^N(p^*) + 1} - \frac{D^{\{i\}}(p)}{D^N(p) + 1}\right), $$

where the price $\phi$ is given by
$$ \phi = \frac{v^{\theta}(N) - v^{\widehat{\theta}}(N)}{ \left( \frac{D^N(p)}{D^N(p)+1} \right)},$$

with $\widehat{\theta} = (N, p, (p_i - 1/\beta)_{i \in N}, \alpha, \beta)$, i.e., $\widehat{\theta}$ is a TC situation with a constant marginal profit for all operators ($1/\beta)$ and with the same total market share as $\theta$ ($D(p)/(D(p)+1)$). Recall from Remark \ref{remark4} that players cannot gain from such a TC situation (because $p_i - c_i = \frac{1}{\beta}$ for all $i,j \in N$). Hence, the numerator of $\phi$ represents the total additional return that is gained compared to a TC situation with the same total market share and where collaborating is not beneficial at all. This total additional return is then divided by the total market share. Indeed, $\phi$ can be recognized as the additional return per unit of market share. 

%This simple price is based on the intuition that the price per dose is equal to the total additional return divided by the total number of available doses. To formalize this intuition, we introduce A(N) as the average additional return per dose when all players in N cooperate.

%\textcolor{black}{From Remark 1, we learned that there exist many TC situations for which the marginal profit is constant. We studied our MSE rule for very marginal profits.}

\begin{table}[h!] \centering
\begin{tabular}{lcccccc} 
            &  \multicolumn{3}{c}{$N$}\\ 
      \cline{2-4}
     & 3 & 4 & 5  \\ \hline
    $MSE$ & 1.0000 & 1.0000 & 1.0000 \\  \hline
\end{tabular}
\caption{Outcome for the Market rule (based on 10,000 TC situations)}
\label{Marketrule}
\end{table}

From Table \ref{Marketrule}, we learn that, for the given TC situations, the allocations of the MSE rule always belong to the core. This does not turn out be a coincidence. Actually, the allocations of the MSE rule are core guaranteed, which will be formalized next.

%\begin{table}[h!] \centering
%\begin{tabular}{lcccccc} 
%      &  \multicolumn{3}{c}{From random prices} &  \multicolumn{3}{c}{From Nash prices}\\ 
%            &  \multicolumn{3}{c}{$N$} &  \multicolumn{3}{c}{$N$}\\ 
%      \cline{2-7}
%     & 3 & 4 & 5 & 3 & 4 & 5  \\ \hline
%    $\Phi(N,v^{\theta})$ & 1.0000 & 1.0000 & 1.0000 & 1.0000 & 1.0000 & 1.0000 \\  \hline
%\end{tabular}
%\caption{Outcome for the Market rule (based on 10,000 TC situations)}
%\label{Marketrule}
%\end{table} 

\begin{theorem} \label{laatstethem} For each TC situation $\theta \in \Theta$ and associated TC game $(N,v^{\theta})$, it holds that $MSE \in \mathscr{C}(N,v^{\theta})$.\end{theorem}

The proof of Theorem \ref{laatstethem} consists of two steps. First, we show the MSE satisfies efficiency, which follows by construction of MSE. Thereafter, we show that MSE satisfies stability. Here we make use of an elementary property of the e-function (see Lemma \ref{lemmaemacht} in Appendix A).

\begin{remark} We could also have formulated other $\widehat{\theta}$'s for which the marginal profit is constant (e.g., we could have selected $\widehat{\theta}$ with $p_i-c_i = \frac{2}{\beta}$ for all $i,j, \in N$). However, for such settings, core non-emptiness cannot be guaranteed.  \end{remark}

We conclude this section with an example, illustrating the calculation of MSE.

\begin{example} Reconsider the TC situation $\theta \in \Theta$ and game $(N,v^{\theta})$ of Example \ref{firstexample1}. Then, $\phi = 3.202$ and 
the allocation of the MSE rule is given by  $MSE = (0.738,	0.296,	0.753)$. One can check that MSE does belong to the core.

\end{example}

\section{An extension of the TC game}\label{extension}

\textcolor{black}{In the introduction, we described that horizontal agreements may qualify for an exemption if they create sufficient pro-consumer benefits that outweigh the anti-competitive effects. Besides, they should not eliminate the competition in the relevant market, implying that participants should have a small market share and their  combined  market  share  should not exceed a specified limit. For some forms of collaboration, it is also necessary to pay back part of the joint profit to society (Article 101(3), TFEU). In this section, we investigate an extended TC game where part of the joint profit can be reallocated to society. In particular, we study which fraction can be reallocated to society, such that the MSE rule, which has proven to be very successful for a stable collaboration, still produces core allocations.} \bigskip

Consider a TC situation $\theta \in \Theta$ and associated game $(N,v^{\theta})$. Now, let $\delta \in (0,1)$ be the fraction of the joint profit ($v^{\theta}(N)$) that players want to reallocate to society. We assume that this applies to any coalition $M \subseteq N$ with $\vert M \vert \geq 2$, i.e., any $M$ with $\vert M \vert \geq 2$ will reallocate $\delta \cdot v^{\theta}(M)$ to society. So, the remaining joint profit of $M$ equals $(1-\delta) \cdot v^{\theta}(M)$. We formalize this new setting in a game $(N,v^{\theta,\delta})$, which we call the TC-$\delta$ game, and it reads as follows

$$v^{\theta,\delta}(M) = \left\{\begin{matrix} (1 - \delta) \cdot v^{\theta}(M) && \mbox{ if }  \vert M \vert \geq 2 \\ \\ v^{\theta}(M) &&\mbox{ if } \vert M \vert \leq 1.\end{matrix} \right.  $$

First of all, observe that there is no reason to collaborate (i.e., the core is empty) if $v^{\theta, \delta}(N) < \sum_{i \in N} v^{\theta, \delta}(\{i\})$. In other words, the fraction that can be paid back at most is

\begin{equation} \label{deltazinvol} \delta \leq 1 - \frac{\sum_{i \in N}v^{\theta}(\{i\})}{v^{\theta}(N)} = 1- \beta \frac{\sum_{i \in N} (p_i - c_i) D^{\{i\}}(N)}{D^N(p^*) \ln \left(  \frac{D^N(c)}{D^N(p)}\right)}. \end{equation}

Naturally, we restrict our attention to TC-$\delta$ games for which (\ref{deltazinvol}) holds true. For these games, we provide a sufficient condition for core non-emptiness. Along with this result, we also provide an allocation rule that produces allocations that belong to the core. \bigskip

\begin{theorem} \label{finalth} If $\delta \leq 1 - \max_{i \in N}\frac{v^{\theta}(\{i\})}{MSE_i}$, then $(1-\delta) MSE \in \mathscr{C}(N,v^{\delta}) \not = \emptyset.$ \end{theorem}

The proof of Theorem \ref{finalth} consists of two steps. First we show that $(1-\delta) MSE$ is efficient, which follows by its construction. Thereafter, we show that $(1-\delta) MSE$ is stable. In doing so, we use that allocations of MSE belong to the core of game $(N,v^{\theta})$.

We conclude this section with an example.

\begin{example}
Reconsider the TC situation of Example \ref{firstexample1} with $\delta = 0.08$. In Table \ref{game11} below, we present the coalitional values for game $(N,v^{\theta,\delta})$.

\begin{table}[h!] \centering
\begin{tabular}{ccccccccc} \hline
    $M$ &  $\{\emptyset\}$ & \{1\} & \{2\} & \{3\} & \{1,2\} & \{1,3\} & \{2,3\} & \{1,2,3\}\\ \hline
    $v^{\theta}(M)$ & 0 & -0.440  &  0.260  &  0.199  &  0.212 &  1.366 &  0.695 & 1.644 \\ \hline
    \end{tabular}
\caption{Coalitional values of game $(N,v^{\theta,\delta})$}
\label{game11}
\end{table}

We have $MSE = (0.738, 0.296,	0.753)$ and so $(1-\delta)MSE=(0.679, 0.272,	0.693)$. It is easy to check that $(1-\delta)MSE$ is a core allocation. We could also have concluded this from Theorem \ref{finalth}, since $\delta = 0.08 < 1 - \max\{\frac{-0.440}{0.738}, \frac{0.260}{0.296}, \frac{0.199}{0.753}\} = 0.124$.\end{example}

\section{Conclusions}
\label{conclu}

In this paper, we introduced a cooperative transport choice (TC) game in which a set of transport operators can collaborate and decide at what prices to offer sustainable urban mobility solutions. To better reflect the decision making process of the travelers, we assume that they chose among the services offered  according to the most widely-used disaggregate demand  model, the multinomial logit model. To be in line with the conditions associated with horizontal agreement exemptions, our TC game assumes that the transport operators optimize their prices, while keeping their total market share constant.

We presented various intuitive allocation rules for our TC game and studied to which extent these allocation rules produce allocations that belong to the core. We showed that two intuitive proportional allocation rules, as well as the well-known Shapley value do not always generate core allocations and therefore cannot sustain the collaboration. We then introduced a market share exchange allocation rule that first allocates to each transport operator the profit he or she generates under collaboration and subsequently compensates those transport operators that lost market share, with additional profit earned by the ones that gained some extra market share. This exchange of market share is facilitated by a unique price, which can be expressed as the additional return by cooperating per unit of market  share. We proved that this allocation rule  sustains the collaboration (i.e., the allocations of the market share exchange rule always belong to the core). Finally, we studied a setting where the transport operators need to pay back part of the joint profit to society.  We showed that, under some natural conditions, the market share exchange rule still sustains the collaboration. We would like to emphasize that most of our results are stable against some deviations in the modelling of players outside a coalition. For instance, if we would use the pessimistic approach of \citet{lardon2019coalitional} (i.e., players outside a coalition select prices that minimize the coalitional profit), the allocations of the market share exchange rule are still core allocations, implying that the core of our TC game is still non-empty.

While the inspiration for our TC game came from the field of urban mobility, results also hold for other applications, for example in city logistics where horizontal collaboration between courier, express and parcel carriers has been recognized as a possible solution to tackle the ‘last-mile’ issue. Accordingly, a first natural direction for future research therefore lies in applying our TC game and its properties to real-world problems (and data). Another direction could be to investigate cooperative games based on more advanced discrete choice models allowing even more complex and precise representations of individual behavior.

%\color{red} \begin{itemize}
%    \item explain that results also hold true for other assumption on coalition: We could also assume that players in $N \backslash M$ set prices that minimize the joint profit of coalition $M$, i.e., the players in $N\backslash M$ would like to sabotage the players in $M$, since they broke up the grand coalition. Under this assumption, core will be larger. As a consequence, numerical study will be a bit more positive (but now 100\%). Theorem 2 and Theorem 3 will still hold.
%\end{itemize}
%  \color{black}

Finally we want to conclude by saying that within the transport community, there is a growing interest in exploiting multidisciplinary methods. By investigating a choice-based cooperative game, we hope that we successfully contributed to bridge the gap between cooperative game theory and discrete choice modelling and that our study can encourage researchers to combine the strengths of these two fields. \bigskip \bigskip 

\noindent \textbf{Acknowledgements} \medskip

The authors would like to thank dr. Ahmadreza Marandi for the fruitful discussion on optimality conditions for non-linear, constraint, optimization problems.

\bibliography{references}

\section*{Appendix A}

In this section, we present the proofs of all theorems and lemmas. For proofs of remarks (Remark 1 and 6) we refer to the document "Supplementary material Appendix A". \bigskip

\noindent \underline{\textbf{Proof of Theorem \ref{optimres}}} \bigskip

This proof consists of three steps. First we identify how our maximization optimization problem relates to another, minimization optimization problem. Then, we identify an optimal price vector and the associated optimal value for this minimization optimization problem. Finally, we relate these outcomes to our original optimization problem.\bigskip

\noindent \underline{\textbf{Step 1. An equivalent optimization problem}} \bigskip

Let $\theta \in \Theta$. First, we replace the constraint of optimization problem $\mathscr{P}$ by another constraint. Recall that the constraint of $\mathscr{P}$ is given by
\begin{equation} \label{dommeconstraint} \frac{\sum_{i \in N} e^{\alpha_i - \beta p_i}}{1+\sum_{j \in N} e^{\alpha_j - \beta p_j}} =  \frac{\sum_{i \in N} e^{\alpha_i - \beta x_i}}{1+\sum_{j \in N} e^{\alpha_j - \beta x_j}}. \end{equation}
Since an e-function is always strictly positive, the above constraint can be rewritten as
$$\left( \sum_{i \in N} e^{\alpha_i - \beta p_i} \right) \left(1 + \sum_{i \in N} e^{\alpha_i - \beta x_i}\right) = \left( \sum_{i \in N} e^{\alpha_i - \beta x_i} \right) \left(1+ \sum_{i \in N} e^{\alpha_i - \beta p_i} \right),$$
which is equal to 
$$ \sum_{i \in N} e^{\alpha_i - \beta p_i}  + \left( \sum_{i \in N} e^{\alpha_i - \beta p_i} \right) \left(\sum_{i \in N} e^{\alpha_i - \beta x_i} \right) =  \sum_{i \in N} e^{\alpha_i - \beta x_i} +  \left( \sum_{i \in N} e^{\alpha_i - \beta x_i} \right) \left( \sum_{i \in N} e^{\alpha_i - \beta p_i} \right).$$
By subtracting $\left( \sum_{i \in N} e^{\alpha_i - \beta p_i} \right) \left(\sum_{i \in N} e^{\alpha_i - \beta x_i} \right)$ on both sides, we obtain:
\begin{equation} \label{handigheidje} \sum_{i \in N} e^{\alpha_i - \beta p_i}  =  \sum_{i \in N} e^{\alpha_i - \beta x_i}, \end{equation}

Hence, we can replace constraint (\ref{dommeconstraint}) by constraint (\ref{handigheidje}).

Now, let $\gamma = 1/\left(1+  \sum_{j \in N} e^{\alpha_j - \beta p_j}\right)$. By using (\ref{handigheidje}), $\mathscr{P}$ can be reformulated as:
\begin{equation} \label{optimizationproblem11}  \begin{aligned}   & \hspace{10mm} \mathscr{P}= \max_{x \in \mathbb{R}^N} \gamma \sum_{i \in N } (x_i-c_i) e^{\alpha_i - \beta x_i} \\& \mbox{ s.t.  }   \sum_{i \in N} e^{\alpha_i - \beta p_i} =  \sum_{i \in N} e^{\alpha_i - \beta x_i}.
\end{aligned} \end{equation}

By using the shorthand notation $D(x) = \sum_{i \in N}e^{\alpha_i - \beta x_i}$  for all $x \in \mathbb{R}$, the above optimization problem can be rewritten as:

\begin{equation*}   \begin{aligned}  & \hspace{10mm} \max_{x \in \mathbb{R}^N} \gamma \sum_{i \in N } (x_i-c_i) e^{\alpha_i - \beta x_i} \\& \mbox{ s.t. }   D(p) =  \sum_{i \in N} e^{\alpha_i - \beta x_i}.
\end{aligned} \end{equation*}
Observe that we could also study the above's optimization problem without constant $\gamma$. In that case, the optimal value would be off a factor $\gamma$ only. Moreover, since maximizing a certain function is the same as minimizing that function times minus one, it is also possible to study the following optimization problem $\mathscr{P}'$, instead of $\mathscr{P}$:

\begin{equation*} \label{optimizationproblem}  \begin{aligned}  & \hspace{10mm} \mathscr{P}' =  \min_{x \in \mathbb{R}^N} -\sum_{i \in N } (x_i-c_i) e^{\alpha_i - \beta x_i} \\& \mbox{ s.t. }   D(p) =  \sum_{i \in N} e^{\alpha_i - \beta x_i}.
\end{aligned} \end{equation*}

Please, note that the optimal value of optimization problem $\mathscr{P}$ equals the optimal value of optimization problem $\mathscr{P}$ times $-\gamma$. Moreover, any optimal price vector of $\mathscr{P}'$ is also an optimal price vector in $\mathscr{P}$. Hence, we can thus also study optimization problem $\mathscr{P}'$. \bigskip

\noindent \underline{\textbf{Step 2. Optimal price vector and associated  optimal value for $\mathscr{P}'$}} \bigskip

In this step, we find the minimal value of $\mathscr{P}'$ and an associated optimal price vector. We do so by applying Theorem \ref{bazaraa}. In terms of Theorem \ref{bazaraa}, we can recognize our optimization problem $\mathscr{P}'$ as a nonlinear optimization problem with objective function
$$f(x) = -\sum_{i \in N } (x_i-c_i) e^{\alpha_i - \beta x_i} \mbox{ for all } x \in \mathbb{R}^N,$$
and equality constraint
$$ h(x) =  \sum_{i \in N} e^{\alpha_i - \beta x_i} - D(p) = 0 \mbox{ for all } x \in \mathbb{R}^N.$$

Moreover, for any $\lambda \in \mathbb{R}$, function $\mathscr{L}(\lambda)$ of Theorem \ref{bazaraa} is given by
\begin{equation*} \mathscr{L}(\lambda) = \min_{x \in \mathbb{R}^N} f(x) + \lambda h(x) = \min_{x \in \mathbb{R}^N} \left[-\sum_{i \in N } (x_i-c_i) e^{\alpha_i - \beta x_i} + \lambda \left(\sum_{i \in N} e^{\alpha_i - \beta x_i} - D(p)\right) \right]. \end{equation*}

In order to apply Theorem \ref{bazaraa}, we first solve $\mathscr{L}(\lambda)$ analytically for any $\lambda \in \mathbb{R}$. Thereafter, we construct a feasible solution $x^* \in \mathbb{R}^N$ and $\lambda^* \in \mathscr{L}$ such that $f(x^*) = \mathscr{L}(\lambda^*$). 
We solve function  $\mathscr{L}(\lambda)$ by dividing the minimization problem in $\vert N \vert$ subproblems. We can do so, since there is no dependency between the variables in $x$. Hence, for any $\lambda \in \mathbb{R}$
\begin{equation} \label{subproblem} \mathscr{L}(\lambda) = \sum_{i \in N } \min_{x_i \in \mathbb{R}} \left[ -(x_i-c_i) e^{\alpha_i - \beta x_i} + \lambda e^{\alpha_i - \beta x_i} \right] - \lambda D(p) \end{equation}

We now solve each subproblem of (\ref{subproblem}). That is, we solve 
$$\min_{x_i \in \mathbb{R}} \left[ -(x_i-c_i) e^{\alpha_i - \beta x_i} + \lambda e^{\alpha_i - \beta x_i} \right] \mbox{ for all } i \in N.$$
We do so by studying the derivative of the objective function of each subproblem. Let $i \in N$. The derivative of the objective function, with respect to $x_i$, is
\begin{equation} \label{derivative}  (\beta (x_i -c_i) - 1 - \lambda \beta ) e^{\alpha_i - \beta x_i} \end{equation}

From \ref{derivative}, we learn that the objective of each subproblem is a decreasing function for $x_i < \frac{1}{\beta} + c_i + \lambda$, constant for $x_i = \frac{1}{\beta} + c_i + \lambda$ and increasing for $x_i> \frac{1}{\beta} + c_i + \lambda$. Hence, the minimum is attained at $x_i = \frac{1}{\beta} + c_i + \lambda$, with objective value
$$  - \frac{1}{\beta}e^{a_ i -  1 - \beta c_i - \beta \lambda}.   $$ By applying the above analysis for each subproblem, we conclude, for each $\lambda \in \mathbb{R}$, that
$$\mathscr{L}(\lambda) = - \frac{1}{\beta} e^{-1 - \lambda \beta }\sum_{i \in N}e^{a_ i - \beta c_i } - \lambda D(p) = - \frac{1}{\beta} e^{-1 - \lambda \beta } D(c) - \lambda D(p).$$

Now, let $\lambda^* = \frac{1}{\beta} \left( \ln\left( \frac{D(c)}{D(p)}\right) - 1\right)$ and $x^*_i = c_i + \frac{1}{\beta} \ln \left( \frac{D(c)}{D(p)}\right)$ for all $i \in N$. By substituting $x^*$ in the objective function of $\mathscr{P}'$, we obtain
\begin{equation} \label{deel1} f(x^*) = - \sum_{i \in N} (x_i^* - c_i) e^{\alpha_i - \beta x_i^*} = -\frac{1}{\beta} \ln \left( \frac{D(c)}{D(p)}\right) \sum_{i \in N}e^{\alpha_i - \beta  c_i + \ln \left( \frac{D(p)}{D(c)}\right)}= - \frac{1}{\beta} \ln \left( \frac{D(c)}{D(p)}\right) D(p). \end{equation}

Moreover, $\mathscr{L}(\lambda^*)$ gives

\begin{equation} \begin{aligned} \label{deel2} \mathscr{L}(\lambda^*) & =- \frac{1}{\beta} e^{-1-\lambda^* \beta } D(c) - \lambda D(p) = - \frac{1}{\beta} \frac{D(p)}{D(c)} D(c) - \lambda^* D(p) = - D(p) \left(\frac{1}{\beta} + \lambda^*\right) \\
& = - D(p)\frac{1}{\beta} \ln \left(\frac{D(c)}{D(p)}\right) \end{aligned} \end{equation}

By combining (\ref{deel1}) and (\ref{deel2}),  we learn that 
$$f(x^*) = \mathscr{L}(\lambda^*) = -  D(p)\frac{1}{\beta} \ln \left(\frac{D(c)}{D(p)}\right).$$ 

Hence, by Theorem \ref{bazaraa}, we can conclude that $x^*$ is an optimal price vector of $\mathscr{P}'$. \bigskip

\noindent \underline{\textbf{Step 3. Back to our original optimization problem}} \bigskip

  In step 1, we learned that $x^*$ is also an optimal price vector of $\mathscr{P}$. Moreover, in step 1, we learned that the optimal value of $\mathscr{P}$ equals the optimal value of $\mathscr{P}'$ times $-\gamma$. Hence, the optimal value of optimization problem $\mathscr{P}$ is
\begin{equation} \label{hulpjevoorlater} -\gamma \cdot -  D(p)\frac{1}{\beta} \ln \left(\frac{D(c)}{D(p)}\right) = \frac{D(p)}{\beta(1+D(p))} \ln\left( \frac{D(c)}{D(p)}\right). \end{equation}
This concludes the proof. $\hfill \square$ \bigskip

\noindent \underline{\textbf{Proof of Theorem \ref{res:lemma}}} \bigskip

Let $\theta \in \Theta$ and $M \subseteq N \backslash \{\emptyset\}$. First, we replace the constraint of the optimization problem of $M$ by another constraint. Recall that the constraint is given by

\begin{equation} \label{hulp3433} \frac{\sum_{i \in M} e^{\alpha_i - \beta p_i}}{1+\sum_{j \in N} e^{\alpha_i - \beta p_i}} = \frac{\sum_{i \in M} e^{\alpha_i - \beta x_i}}{1+ \sum_{ i \in M} e^{\alpha_i - \beta x_i} + \sum_{ j \in N \backslash M} e^{\alpha_i - \beta p_i}}. \end{equation}

Since an e-function is always strictly positive, we also have
$$ \left(\sum_{i \in M} e^{\alpha_i - \beta p_i}\right) \left( 1+ \sum_{ i \in M} e^{\alpha_i - \beta x_i} + \sum_{ j \in N \backslash M} e^{\alpha_i - \beta p_i}\right)  = \left( \sum_{i \in M} e^{\alpha_i - \beta x_i} \right)     \left( 1+\sum_{j \in N} e^{\alpha_i - \beta p_i} \right) .$$
By subtracting $(\sum_{i \in M} e^{\alpha_i - \beta p_i})(\sum_{i \in M}e^{\alpha_i - \beta x_i})$ on both sides, we obtain
$$ \left(\sum_{i \in M} e^{\alpha_i - \beta p_i}\right) \left( 1+  \sum_{ j \in N \backslash M} e^{\alpha_i - \beta p_i}\right)  = \left(\sum_{i \in M} e^{\alpha_i - \beta x_i}\right) \left( 1+  \sum_{ j \in N \backslash M} e^{\alpha_i - \beta p_i}\right) .$$
Observe that $(1+ \sum_{i \in N \backslash M} e^{\alpha_i - \beta p_i})$ is a strictly positive constant, and so we can divide both sides of the last equation by this term. This leads to

\begin{equation} \label{hulp22}  \sum_{i \in M} e^{\alpha_i - \beta p_i}  = \sum_{i \in M} e^{\alpha_i - \beta x_i}. \end{equation}

Hence, we can replace constraint (\ref{hulp3433}) by constraint (\ref{hulp22}). \bigskip

\noindent Now, let $\gamma = 1/(1 + \sum_{j \in N} e^{\alpha_i - \beta p_j})>0$. By (\ref{hulp22}), our optimization problem becomes

\begin{equation*} \label{optimizationproblem}  \begin{aligned}   & \hspace{10mm}  \max_{x \in \mathbb{R}^M} \gamma \sum_{i \in M } (x_i-c_i) e^{\alpha_i - \beta x_i} \\& \mbox{ s.t.  }   \sum_{i \in M} e^{\alpha_i - \beta p_i} =  \sum_{i \in M} e^{\alpha_i - \beta x_i}.
\end{aligned} \end{equation*}

Above optimization problem is exactly equal to optimization problem (\ref{optimizationproblem11}) of the proof of Theorem \ref{optimres}, except that we consider set $M$ in stead of $N$. So, we can now use Theorem \ref{optimres} (see equation (\ref{hulpjevoorlater})), to conclude that the optimal value of our optimization problem equals
$$ -\gamma \cdot -D^M(p) \frac{1}{\beta} \ln \left( \frac{D^M(c)}{D^M(p)}\right) = \frac{D^M(p)}{\beta(D^N(p) + 1)} \ln \left( \frac{D^M(c)}{D^M(p)}\right).$$
Hence, $v^{\theta}(M) = \frac{D^M(p)}{\beta(D^N(p) + 1)} \ln \left( \frac{D^M(c)}{D^M(p)}\right)$, which concludes the proof. 

%$$e^{A-B} \geq \left(\frac{A}{B}\right)^B$$
\bigskip

 \begin{lemma} \label{lemmaemacht} For all $A,B \in \mathbb{R}_{++}$, it holds that
$$A \geq B \left(  \ln \left( \frac{A}{B} \right) +1 \right)$$
\end{lemma}

\noindent \textbf{Proof :}  First, we prove that $e^y - ye \geq 0$ for all $y \in \mathbb{R}$. We do so by showing that the function is convex, and identifying that the minimal value of this function equals 0. The function is convex, since $\frac{d^2}{d^2y}(e^y - ye) = e^y \geq 0$. Moreover, we have  $\frac{d}{dy}(e^y - ye) = e^y - e$, implying that the minimal value is attained at $y=1$ with associated function value $e^1-1 \cdot e = 0$. \bigskip

Now, let $A,B \in \mathbb{R}_{++}$. Note that $\frac{A}{B}$ exists, because $A,B >0$. We just learned that $e^y - ye \geq 0$ for all $y \in \mathbb{R}$, and so, we also have
\begin{equation*}  \begin{aligned}  e^\frac{A}{B} &\geq \frac{A}{B}e \\
\Longleftrightarrow{} (e^{\frac{A}{B}})^B &\geq \left(\frac{A}{B}e\right)^{B} \\
%\Longleftrightarrow{}(e^{\frac{A}{B}})^B e^{-B} &\geq \left(\frac{A}{B}\right)^{B} \\
%\Longleftrightarrow{}e^{A-B} &\geq \left(\frac{A}{B}\right)^{B} \\
%\Longleftrightarrow{} e^A &\geq e^{ \ln \left( \left( \frac{eA}{B}\right)^B \right)} \\
\Longleftrightarrow{} e^A &\geq e^{ \ln \left( \left(\frac{A}{B}e \right)^B\right)} \\
\Longleftrightarrow{} A &\geq B  \ln \left( \frac{eA}{B}\right) \\
\Longleftrightarrow{} A &\geq B  \left( \ln \left( \frac{A}{B}\right) + 1  \right) 
\end{aligned}  \end{equation*}

where the first implication holds since $e^{\frac{A}{B}}>0$, $\frac{A}{B}e >0$ and $B>0$. The second implication follows from the fact that $e^{\ln{x}}=x$ for $x \in \mathbb{R}_{++}$ and $(e^{x})^y = e^{xy}$ for all $x,y \in \mathbb{R}$. The third implication results from the fact that $e>0$ and $\ln(x^y) = y \ln (x)$ for all $x,y \in \mathbb{R}_{++}$. The last implication is a result of property $\ln(xy) = \ln(x) + \ln(y)$ for all $x,y \in \mathbb{R}_{++}$ and the fact that $\ln(e)=1$. This concludes the proof.
$\hfill \square$ \bigskip

\noindent \underline{\textbf{Proof of Theorem \ref{laatstethem}}} \bigskip

Let $\theta \in \Theta$ and consider the associated game $(N,v^{\theta})$. We show that allocation $MSE \in \mathscr{C}(N,v^{\theta})$. We do so by showing that

$$\sum_{i \in N} MSE_i = v^{\theta}(N),$$
$$\sum_{i \in M} MSE_i \geq v^{\theta}(M) \mbox{ for all } M \subseteq N.$$

For the first part, recall that
$$MSE_i = (p^*-c_i) \frac{D^{\{i\}}(p^*)}{D^N(p^*)+1} - \phi \left(\frac{D^{\{i\}}(p^*)}{D^N(p^*)+1} - \frac{D^{\{i\}}(p)}{D^N(p)+1}\right) \mbox{ for all } i \in N.$$

Now, observe that
\begin{equation*} \begin{aligned} \sum_{i \in N} MSE_i &= \sum_{i \in N }   \left( (p^*_i-c_i) \frac{D^{\{i\}}(p^*)}{D^{N}(p^*)+1} - \phi \left(\frac{D^{\{i\}}(p^*)}{D^N(p^*)+1} - \frac{D^{\{i\}}(p)}{D^N(p)+1}\right) \right) \\
&= \sum_{i \in N}(p^*_i-c_i) \frac{D^{\{i\}}(p^*)}{D^N(p^*)+1} - \phi \left( \sum_{i \in N} \frac{D^{\{i\}}(p^*)}{D^N(p^*)+1} - \sum_{i \in N} \frac{D^{\{i\}}(p)}{D^N(p)+1} \right)\\
& = \frac{1}{\beta} \ln \left( \frac{D^N(c)}{D^N(p)}\right) \frac{D^N(p)}{D^{N}(p)+1} \\
&= v^{\theta}(N). \end{aligned} \end{equation*}

In the third equality, we use the definition of $p^*$, apply that $\sum_{i \in N}D^{\{i\}}(p^*) = D^N(p^*)$ and consequently use that $\frac{D^N(p^*)}{D^N(p^*) + 1} = \frac{D^N(p)}{D^N(p)+1}$ (i.e., the total market share remains stable). \bigskip

\noindent Now, it remains to prove $\sum_{i \in M} MSE_i \geq v^{\theta}(M) \mbox{ for all } M \subseteq N$. First, observe that,

\begin{equation} \label{phiomschrijven} \phi = \frac{v^{\theta}(N) - v^{\theta^*}(N)}{ \left( \frac{D^N(p)}{D^N(p)+1} \right)} = \frac{1}{\beta} \ln\left( \frac{D^N(c)}{D^N(p)}\right) - \frac{1}{\beta}. \end{equation}

Next, by exploiting $p^*$, we learn that
\begin{equation} \label{Dngelijkaandnp} D^N(p^*) = \sum_{i \in N} e^{\alpha_i - \beta p^*_i} = \sum_{i \in N} e^{\alpha_i - \beta \left(c_i + \frac{1}{\beta} \ln \left( \frac{D^N(c)}{D^N(p)}\right)\right)} = D^N(c) \cdot \frac{D^N(p)}{D^N(c)} = D^N(p).  \end{equation}

By using (\ref{phiomschrijven}), (\ref{Dngelijkaandnp}) and exploiting $p^*$, we can reformulate $MSE_i$ for all $i \in N$ as follows
\begin{equation}  \begin{aligned} \label{hulpje0}
MSE_i &= (p^*_i-c_i)  \frac{D^{\{i\}}(p^*)}{D^N(p^*)+1} - \phi \left(\frac{D^{\{i\}}(p^*)}{D^N(p^*) + 1} - \frac{D^{\{i\}}(p)}{D^N(p) + 1}\right) \\
&= \frac{1}{D^N(p) + 1} \bigg[(p^*_i-c_i) D^{\{i\}}(p^*) - \phi \left(D^{\{i\}}(p^*) - D^{\{i\}}(p)\right) \bigg] \\
&= \frac{1}{D^N(p) + 1} \bigg[\frac{1}{\beta }\ln\left(\frac{D^N(c)}{D^N(p)} \right) \frac{D^N(p)}{D^N(c)} D^{\{i\}}(c)  
\\ & \hspace{26mm} - \left( \frac{1}{\beta} \ln\left( \frac{D^N(c)}{D^N(p)}\right) - \frac{1}{\beta} \right) \left(\frac{D^N(p)}{D^N(c)} D^{\{i\}}(c) - D^{\{i\}}(p)\right) \bigg] \\
& = \frac{1}{D^N(p) + 1} \bigg[ \frac{1}{\beta} \ln\left( \frac{D^N(c)}{D^N(p)}\right)  D^{\{i\}}(p) + \frac{1}{\beta} \left( \frac{D^N(p)}{D^N(c)} D^{\{i\}}(c) - D^{\{i\}}(p)\right) \bigg] \\
& = \frac{1}{\beta(D^N(p) + 1)} \bigg[  \ln\left( \frac{D^N(c)}{D^N(p)}\right)  D^{\{i\}}(p) + \frac{D^N(p)}{D^N(c)} D^{\{i\}}(c) - D^{\{i\}}(p) \bigg]. \end{aligned}
\end{equation}

Moreover, from Lemma \ref{lemmaemacht}, we learned that
\begin{equation} \label{hulpje1} A \geq B  \left( \ln \left( \frac{A}{B}\right) + 1  \right)
\end{equation}

Now, fix an $M \subseteq N$. Moreover, let's use the inequality of (\ref{hulpje1}), and set $A = \frac{D^M(c)D^N(p)}{D^N(c)}$ and $B = D^M(p) >0$. Then,
\begin{equation} \label{hulpje2} \begin{aligned}
&\frac{ D^M(c) D^N(p)}{D^N(c)} \geq D^M(p) \left( \ln \left( \frac{D^M(c)D^N(p)}{D^N(c) D^M(p)}  \right) + 1\right)  \\
&\Longleftrightarrow{}  \frac{ D^M(c) D^N(p)}{D^N(c)} \geq D^M(p) \left( \ln \left( \frac{D^M(c)}{D^M(p)}\right) + \ln \left( \frac{D^N(p)}{D^N(c)}\right) + 1\right)  \\
&\Longleftrightarrow{} 
D^M(p) \ln \left( \frac{D^N(c)}{D^N(p)}\right) + \frac{ D^M(c) D^N(p)}{D^N(c)} - D^M(p) \geq D^M(p) \ln \left( \frac{D^M(c)}{D^M(p)}\right) \\ 
&\Longleftrightarrow{} 
\frac{1}{\beta (D^N(p)+1)} \left( D^M(p) \ln \left( \frac{D^N(c)}{D^N(p)}\right) + \frac{D^N(p)}{D^N(c)} D^M(c)  - D^M(p) \right) \\
& \hspace{10mm} \geq \frac{D^M(p)}{\beta (D^N(p)+1)} \ln \left( \frac{D^M(c)}{D^M(p)}\right) \\ 
&\Longleftrightarrow{} 
\frac{1}{\beta  (D^N(p)+1)} \left( \ln \left( \frac{D^N(c)}{D^N(p)}\right) D^M(p) + \frac{D^N(p)}{D^N(c)} D^M(c)  - D^M(p) \right) \geq v^{\theta}(M)
\end{aligned}
\end{equation}

Please, note that in the last inequality we used the definition of $v^{\theta}(M)$. By using the last equality of (\ref{hulpje0}) and the last inequality of (\ref{hulpje2}), we have
\begin{equation*}
    \begin{aligned}
    \sum_{i \in M} MSE_i &= \sum_{i \in M} \frac{1}{\beta (D^N(p)+1)}\left(\ln\left(\frac{D^N(c)}{D^N(p)} \right)D^{\{i\}}(p) +   \frac{D^N(p)}{D^N(c)}D^{\{i\}}(c) -D^{\{i\}}(p) \right) \\
&=  \frac{1}{\beta (D^N(p)  +1)}\left(\ln\left(\frac{D^N(c)}{D^N(p)} \right) \sum_{i \in M}D^{\{i\}}(p)  +  \frac{D^N(p)}{D^N(c)} \sum_{i \in M}D^{\{i\}}(c)  -\sum_{i \in M} D^{\{i\}}(p) \right) \\
&=  \frac{1}{\beta (D^N(p)  +1)}\left(\ln\left(\frac{D^N(c)}{D^N(p)} \right)D^M(p) +  \frac{D^N(p)}{D^N(c)}D^M(c)  -D^M(p) \right) \\
& \geq v^{\theta}(M),
\end{aligned} \end{equation*}
which is exactly what we need to show. This concludes the proof. $\hfill \square$ \bigskip

\noindent \underline{\textbf{Proof of Theorem \ref{finalth}}} \bigskip

Let $\delta \leq 1 - \max_{i \in N} \frac{v^{\theta}(\{i\})}{MSE_i}$. We will show that allocation rule $(1-\delta)MSE \in \mathscr{C}(N, v^{\theta, \delta}) \not = \emptyset$. We do so by showing that $(1-\delta)MSE$ satisfies

$$\sum_{i \in N} (1-\delta)MSE_i = v^{\theta, \delta}(N),$$
$$\sum_{i \in M} (1-\delta) MSE_i \geq v^{\theta, \delta}(M) \mbox{ for all } M \subseteq N.$$

For the first part, observe that
$$ \sum_{i \in N} (1-\delta)MSE_i =  (1- \delta) \sum_{i \in N} MSE_i = (1- \delta) v^{\theta}(N) = v^{\theta, \delta}(N),$$

where the second equality holds since $MSE \in \mathscr{C}(N,v^{\theta})$ (See Theorem \ref{laatstethem}).

Now, it remains to prove that $\sum_{i \in M} (1-\delta) MSE_i \geq v^{\theta, \delta}(M) \mbox{ for all } M \subseteq N$. Since $\delta \leq 1 - \max_{i \in N} \frac{v^{\theta}(\{i\})}{MSE_i}$, we also have $\delta \leq 1 - \frac{v^{\theta}(\{i\})}{MSE_i}$ for all $i \in N$. From this, we can conclude that $(1-\delta)MSE_i \geq v^{\theta}(\{i\}) = v^{\theta,\delta}(\{i\})$ for all $i \in N$. Hence, it remains to show that $\sum_{i \in M}(1-\delta)MSE_i \geq v^{\theta, \delta}(M)$ for all $M \subseteq N$ with $ \vert M \vert \geq 2$. 

Let $M \subseteq N$ with $ \vert M \vert \geq 2$. We have

$$ \sum_{i \in M}(1-\delta)MSE_i = (1-\delta) \sum_{i \in M }MSE_i \geq (1-\delta) v^{\theta}(M) = v^{\theta, \delta}(M),$$

where the inequality holds since $MSE \in \mathscr{C}(N,v^{\theta})$ (see Theorem \ref{laatstethem}). This implies that $\sum_{i \in M} MSE^{\delta}_i \geq v^{\theta, \delta}(M)$ for all $M \subseteq N$. This concludes the proof. $\hfill \square$

\begin{theorem} \label{bazaraa}
Let $f: \mathbb{R}^N \to \mathbb{R}$ and
$h(x): \mathbb{R}^N \to \mathbb{R}$. If $x^* \in \mathbb{R}^N$ is a feasible solution of 
nonlinear programming problem $\mathscr{G}$:
\begin{equation*} \begin{aligned} \mathscr{G} = \min f(x) \\
h(x) = 0\\
x \in \mathbb{R}^N,
\end{aligned} \end{equation*}
and $f(x^*) = \mathscr{L}(\lambda)$ for some $\lambda \in \mathbb{R}$, where 
$$\mathscr{L}(\lambda) = \min_{x \in \mathbb{R}^N} \left\{f(x) + \lambda h(x)\right\}$$ 
then $x^*$ is an optimal solution of $\mathscr{G}$.
\end{theorem}

\noindent \textbf{Proof:} See \citet{bazaraa2013nonlinear}, chapter 6, corollary 2. $\hfill \square$ \bigskip

\newpage

\section*{Supplementary material Appendix A}

\noindent \underline{\textbf{Proof of Remark 1}} \bigskip

For each TC situation $\theta \in \Theta$ for which $p$ is a Nash equilibrium, we show that
\begin{equation}
p_i = \frac{1+W\left(\frac{e^{\alpha_i-1-\beta c_i}}{1+ \sum_{j \neq i} e^{\alpha_j - p_j\beta} }\right)}{\beta} + c_i \hspace{2mm} \mbox{ for all } i \in N. \label{watwewillen1}
\end{equation} 

Let $\theta \in \Theta$ and $p$ be a Nash equilibrium. Hence, $p$ is also a solution of the set of first-order conditions, based on the profit functions of the individual transport operators. We will now derive these first-order conditions and show that they coincide with equation (\ref{watwewillen1}).

The derivative of the profit function of transport operator $i \in N$, with respect to price (which we denote by $p_i'$ instead of $p_i$ for notational convenience), equals
%We will do so by showing that $p$ is a solution of the set of first-order conditions, based on the profit functions of the individual transport operators. Let $\theta \in \Theta$ and $i \in N$. Then, the derivative of the profit function of transport operator $i$ with respect to $p_i$ is

\begin{equation*} \begin{aligned} &\frac{d}{dp_i'} \left( (p'_i - c_i) \cdot \left( \frac{e^{\alpha_i - \beta p_i'}}{1+\sum_{j \in N}e^{\alpha_j - \beta p'_j}} \right) \right)  \\
&= \frac{e^{\alpha_i - \beta p'_i}}{1+\sum_{j \in N}e^{\alpha_j - \beta p'_j}} \left( 1 - \beta \left( \frac{1 + \sum_{j \in N \backslash \{i\}}e^{\alpha_j - \beta p'_j}}{1+\sum_{j \in N}e^{\alpha_j - \beta p'_j}}\right) (p'_i - c_i)\right). \end{aligned} \end{equation*}

Since $\frac{e^{\alpha_k - \beta p'_k}}{1+\sum_{j \in N}e^{\alpha_j - \beta p'_j}} > 0$ for all $k \in N$, the first-order conditions (with Nash equilibrium prices $p$) read as follows

$$  1 - \beta \left(1 + \frac{\sum_{j \in N \backslash \{i\}}e^{\alpha_j - \beta p_j}}{1+\sum_{j \in N}e^{\alpha_j - \beta p_j}}\right) (p_i - c_i) = 0 \mbox{ for all } i \in N. $$

These first-order conditions can be written in terms of the Nash equilibrium prices $p$

% \textbf{Proof:} The Nash prices are obtained by solving:
%\begin{equation*}
%\max_{p_i \in \mathbb{R}} (p_i - c_i) \cdot S_i
%\end{equation*} 
%where $S_i = \left( \frac{e^{\alpha_i - \beta p_i}}{1+\sum_{j \in N}e^{\alpha_j - \beta p_j}} \right)$ denotes the market share of operator $i$. 

%\noindent From first-order conditions, we get:
%\begin{align}
%    p_i &= \frac{1}{\beta \cdot (1-S_i)} +c_i \label{eq1}
%\end{align}

\begin{equation}   p_i = \frac{1}{\beta \cdot \frac{1 + \sum_{j \in N \backslash \{i\}}e^{\alpha_j - \beta p_j}}{1+\sum_{j \in N}e^{\alpha_j - \beta p_j}}}   +c_i \mbox{ for all } i \in N. \label{eq1}\end{equation}

We continue by rewriting equation (\ref{eq1}) towards equation (\ref{watwewillen1}). For all  $i \in N$,  we have
%\noindent Replacing market shares by their expressions, we obtain:
\begin{align}
    p_i &= \frac{1}{\beta \cdot \left(1 - \frac{e^{\alpha_i - \beta p_i}}{1+\sum_{j \in N}e^{\alpha_j - \beta p_j}} \right)}   +c_i \\
  %  p_i &= \frac{1}{\beta - \beta \cdot \left( \frac{e^{\alpha_i - \beta p_i}}{1+\sum_{j \in N}e^{\alpha_j - \beta p_j}} \right)} +c_i  \\
  &= \frac{1}{\frac{\beta \cdot (1+\sum_{j \in N}e^{\alpha_j - \beta p_j} ) - \beta \cdot (e^{\alpha_i-\beta p_i})}{1+\sum_{j \in N}e^{\alpha_j - \beta p_j} }} +c_i \\ 
      &= \frac{1+\sum_{j \in N}e^{\alpha_j - \beta p_j}}{\beta \cdot (1+ \sum_{j \in N \backslash \{ i\}} e^{\alpha_j - \beta p_j}  )} + c_i
    \end{align}
  
  \begin{align}
    &= \frac{1+ \sum_{j \in N \backslash \{i\}} e^{\alpha_j - \beta p_j} + e^{\alpha_i - \beta p_i}}{\beta \cdot (1+ \sum_{j \in N \backslash \{i\}} e^{\alpha_j - \beta p_j}  )} + c_i\\
     %    &=  \frac{1}{\beta} +  \frac{e^{\alpha_i - \beta p_i}}{\beta \cdot (1+ \sum_{j \in N \backslash \{i\}} e^{\alpha_j - \beta p_j}  )} + c_i\\       
    &= \frac{1}{\beta} + \frac{e^{\alpha_i - \beta p_i}}{\beta \cdot A} +c_i \label{single1}
\end{align}

\noindent with $A = 1 + \sum_{j \in N \backslash \{ i\}} e^{\alpha_j - \beta p_j}$. \\

\noindent Multiplying equation (\ref{single1}) by $\beta$ and then subtracting $\alpha_i$, we obtain
\begin{align}
  \beta p_i - \alpha_i  &= 1 +  \frac{e^{\alpha_i - \beta p_i}}{A} + \beta c_i - \alpha_i \\
\Longleftrightarrow{} \frac{e^{\alpha_i - \beta p_i}}{A} - \beta p_i + \alpha_i  &= \alpha_i -1 - \beta c_i \label{single2} 
\end{align}
Taking exponential on both sides of equation (\ref{single2}) and using that $A>0$, we have
\begin{align}
  e^{\frac{e^{\alpha_i - \beta p_i}}{A} - \beta p_i + \alpha_i}  &= e^{\alpha_i -1 - \beta c_i} \\
  \Longleftrightarrow{}  e^{\frac{e^{\alpha_i - \beta p_i}}{A}} \cdot e^{\alpha_i - \beta p_i}  &= e^{\alpha_i -1 - \beta c_i} \\
 \Longleftrightarrow{}   e^{\frac{e^{\alpha_i - \beta p_i}}{A}} \cdot \frac{e^{\alpha_i - \beta p_i}}{A}  &= \frac{e^{\alpha_i -1 - \beta c_i}}{A}, \label{single3}
\end{align}

Now, observe that equation (\ref{single3}) can reformulated in terms of the classic LambertW equation (i..e, as $e^W \cdot W = c$ for some $c \in \mathbb{R}$). In particular,
(\ref{single3}) can be reformulated as

\begin{equation}
  \frac{e^{\alpha_i - \beta p_i}}{A} = W\left(\frac{e^{\alpha_i -1 -  \beta c_i }}{A}\right) \label{single4}
\end{equation}
Taking logarithms on both side of (\ref{single4}), which is allowed since $\frac{e^{\alpha_i - \beta p_i}}{A} >0$, we have
\begin{equation}
  \ln(\frac{e^{\alpha_i - \beta p_i}}{A}) = \ln\left(W\left(\frac{e^{\alpha_i -1 - \beta c_i}}{A}\right)\right) \label{single5}
\end{equation}
Using the logarithmic property of the LambertW function (i.e., $\ln(W(x)) = \ln(x) - W(x)$ for any $x \in \mathbb{R}_{++}$) and the fact that $e^{\alpha - 1 - \beta c_i}/A, e^{\alpha_ i - \beta p_i}/A >0$, equation (\ref{single5}) becomes:
\begin{align}
\ln\left(\frac{e^{\alpha_i - \beta p_i}}{A}\right) &= \ln\left(\frac{e^{\alpha_i -1 -  \beta c_i}}{A}\right) - W\left(\frac{e^{\alpha_i -1 -  \beta c_i }}{A}\right) \\
\Longleftrightarrow{} \alpha_i - \beta p_i -\ln(A) &= \alpha_i -1 - \beta c_i - \ln(A) - W\left(\frac{e^{\alpha_i -1 -  \beta c_i}}{A}\right) \\
\Longleftrightarrow{} - \beta p_i  &=  -1  - \beta c_i - W\left(\frac{e^{\alpha_i -1 - \beta c_i}}{A}\right) \end{align}

\begin{align}
\Longleftrightarrow{} \beta p_i  &= 1 + \beta c_i + W\left(\frac{e^{\alpha_i -1 - \beta c_i}}{A}\right) \\
\Longleftrightarrow{} p_i  &= \frac{1 + W\left(\frac{e^{\alpha_i -1 - \beta c_i}}{A}\right)}{\beta} + c_i
\end{align}
Substituting $A = 1 + \sum_{j \neq i} e^{\alpha_j - \beta p_j}$, we obtain
\begin{equation}
    p_i  = \frac{1 + W(\frac{e^{\alpha_i -1 - \beta c_i}}{1 + \sum_{j \neq i} e^{\alpha_j - \beta p_j}})}{\beta} + c_i,
\end{equation}
which we needed to show. This concludes the proof. $\hfill \square$ \bigskip

\noindent \textbf{\underline{Proof of Remark \ref{delaatsteremarkofzo}}} \bigskip

First we will show that our game is superadditive, i.e., $v^{\theta}(M) + v^{\theta}(K) \leq v^{\theta}(M \cup K)$ for all $M, K \subseteq N$ with $ M \cap K = \emptyset$ and all $\theta \in \Theta$. Let $\theta \in \Theta$ and $M, K \subseteq N$ with $M \cap K = \emptyset$. Moreover, let $x^M = (x^M_i)_{i \in M}$ be an optimal solution of the optimization problem of coalition $M$. Similarly, let $x^K = (x^K_i)_{i \in K}$ be an optimal solution of the optimization problem of coalition $K$. Next, let

$$x^{M \cup K}_i = \left\{ \begin{matrix} x^M_i && \mbox{ if } i \in M \\ \\ x^K_i && \mbox{ if } i \in K. \end{matrix} \right.$$

We will show that $x^{M \cup K} = (x^{M \cup K}_i)_{i \in M \cup K}$ is a feasible solution of the optimization problem of coalition $M \cup K$. That means, we need to show that the market share constraint of the optimization problem of coalition $M \cup K$ is satisfied. Recall that, due to equation (\ref{hulp22}) of Theorem \ref{res:lemma}, the total market share constraint of any coalition $T \subseteq N$ reads as

\begin{equation} \label{handyhandy33} \sum_{i \in T} e^{\alpha_ i - \beta x_i} = \sum_{i \in T} e^{\alpha_ i - \beta p_i}. \end{equation}

Hence, for optimal solutions $x^M$ and $x^K$, we have

\begin{equation} \label{handyhdandy2} \begin{aligned} \sum_{i \in M} e^{\alpha_ i - \beta x_i^M} = \sum_{i \in M} e^{\alpha_ i - \beta p_i} \\
\sum_{i \in K} e^{\alpha_ i - \beta x_i^K} = \sum_{i \in K} e^{\alpha_ i - \beta p_i}. \end{aligned} \end{equation}

As a consequence, we have
$$ \sum_{i \in M \cup K} e^{\alpha_ i - \beta x_i^{M \cup K}} = \sum_{i \in M } e^{\alpha_ i - \beta x_i^{M}} +  \sum_{i \in K} e^{\alpha_ i - \beta x_i^{K}} = \sum_{i \in M} e^{\alpha_ i - \beta p_i} + \sum_{i \in K} e^{\alpha_ i - \beta p_i} = \sum_{i \in M \cup K} e^{\alpha_ i - \beta p_i}. $$

From the above equation, we learn that solution $x^{M \cup K}$ is a feasible solution of the optimization problem of coalition $M \cup K$. Next, we will show that the sum of the objective functions of the optimization problems of $M$ and $K$, evaluated for $x^M$ and $x^K$ coincide with the objective function of the optimization problem of $M \cup K$, evaluated at $x^{M \cup K}$.

The sum of the objective functions of the optimization problems of $M$ and $K$, evaluated for $x^M$ and $x^K$ reads as follows

\begin{equation*} \begin{aligned} & \frac{ \sum_{i \in M} (x^M_i - c_i) e^{\alpha_i - \beta x_i^M}}{ 1 + \sum_{j\in M} e^{\alpha_j - \beta x_j^M} + \sum_{j \in N \backslash M} e^{\alpha_j - \beta p_j}} +  \frac{\sum_{i \in K} (x^K_i - c_i) e^{\alpha_i - \beta x_i^K}}{ 1 + \sum_{j\in K} e^{\alpha_j - \beta x_j^K} + \sum_{j \in N \backslash K} e^{\alpha_j - \beta p_j}}  \\ 
=&  \frac{ \sum_{i \in M} (x^M_i - c_i) e^{\alpha_i - \beta x_i^M}}{ 1 + \sum_{j\in N} e^{\alpha_j - \beta p_j}} + \frac{ \sum_{i \in K} (x^K_i - c_i) e^{\alpha_i - \beta x_i^K}}{ 1 + \sum_{j\in N} e^{\alpha_j - \beta p_j} } \\ 
=& \frac{ \sum_{i \in M \cup K} (x^{M \cup K}_i - c_i) e^{\alpha_i - \beta x_i^{M\cup K}}}{ 1 + \sum_{j\in N} e^{\alpha_j - \beta p_j}} \\
=& \frac{ \sum_{i \in M \cup K} (x^{M \cup K}_i - c_i) e^{\alpha_i - \beta x_i^{M \cup K}}}{ 1 + \sum_{j\in M \cup K} e^{\alpha_j - \beta x_j^{M \cup K}} +  \sum_{j\in N \backslash (M \cup K)} e^{\alpha_j - \beta p_j}}
\end{aligned} \end{equation*}

 Note that we used equation (\ref{handyhandy33}) for coalitions $M$ and $K$ in the first equality and for coalition $M \cup K$ in the last equality. Next, observe that the expression in the last equation coincides with the objective function of optimization problem $M \cup K$, evaluated for $X^{M \cup K}$.

Hence, the sum of the coalitional values $v^{\theta}(M)$ and $v^{\theta}(K)$ coincides with the objective value of the optimization problem of $M \cup K$, evaluated as $X^{M \cup K}$. Since solution $x^{M \cup K}$ is as feasible solution, we conclude that $v^{\theta}(M \cup K)$ is at least this value. Hence,

$$v^{\theta}(M) + v^{\theta}(K) \leq v^{\theta}(M \cup K), $$

which concludes the proof for superadditivty. \bigskip

Now we will provide an example that illustrates that our TC game is not convex nor monotonic in general. Consider a TC situation $\theta \in \Theta$ with $p = (0.5, 0.5, 2)$, $c= (0.5, 1, 1.5)$, $\alpha  = (1,2,1.5)$, and $\beta = 0.1$. The coalitional values are represented in Table \ref{game10} below.

\begin{table}[h!] \centering
\begin{tabular}{ccccccccc} \hline
    $M$ &  $\{\emptyset\}$ & \{1\} & \{2\} & \{3\} & \{1,2\} & \{1,3\} & \{2,3\} & \{1,2,3\}\\ \hline
    $v^{\theta}(M)$ & 0 & 0   &  -0.246  &  0.128  & -0.244   &  0.130 &  -0.109  & -0.109 \\ \hline
    \end{tabular}
\caption{Coalitional values of game $(N,v^{\theta})$}
\label{game10}
\end{table} 

Observe that $v^{\theta}(\{1\}) = 0 > -0.244 = v^{\theta}(\{1,2\})$, implying that the game is not monotonic. In addition, observe that $v^{\theta}(\{1,2\}) - v^{\theta}(\{2\})  = 0.002 > 0 = v^{\theta}(\{1,2,3\}) - v^{\theta}(\{2,3\})$, implying that the game is also not convex.

\end{document}